\documentclass[USenglish,onecolumn]{article}

\usepackage[small]{dgruyter}
\usepackage{amssymb,amsmath,color,soul,enumerate}

\newcommand{\F}{\mathcal{F}}
\newcommand{\RR}{\mathbb{R}}
\newcommand{\D}{\mathcal{D}}
\newcommand{\HH}{\mathcal{H}}
\newcommand{\dist}{\operatorname{dist}}
\newcommand{\lF}{\lambda_{\F}}
\newcommand{\slope}[1]{\lambda_{\F}(#1)}

\numberwithin{equation}{section}

\theoremstyle{plain}
\newtheorem{theorem}{Theorem}[section]

\newtheorem{lemma}[theorem]{Lemma}
\newtheorem{corollary}[theorem]{Corollary}

\theoremstyle{definition}
\newtheorem{definition}[theorem]{Definition}
\newtheorem{example}[theorem]{Example}
\newtheorem{remark}[theorem]{Remark}

\begin{document}

\articletype{Research Article}

\author[1]{Cristian Bereanu}
\author[2]{Nicu\c sor Costea}

\affil[1]{Institute of Mathematics "Simion Stoilow" of the Romanian Academy,
21 Calea Grivi\c tei Street, 010702 Bucharest, Romania.

and

Faculty of Mathematics, University of Bucharest, 14 Academiei Street, 70109 Bucharest, Romania.

E-mail: cbereanu@imar.ro}

\affil[2]{Department of Mathematics and Computer Science,
National University of Science and Technology POLITEHNICA Bucharest,
   313 Splaiul Independen\c tei, 060042 Bucharest, Romania.

E-mail: nicusor.costea2606@upb.ro}

  \title{\large Double interaction principle for the Lorentz force equation}
  \runningtitle{Double interaction principle for the Lorentz force equation}
 \begin{abstract}
{In this paper we study the Lorentz force equation using the Poincar\'e action functional. We show that,
if the scalar electric potential has 'good' interactions  at infinity and around the origin with the vector 
magnetic potential, then the Lorentz force equation has at least two periodic solutions with a fixed period.
For the proof of our main result we introduce a new minimax theorem for the Poincar\'e action functional 
that does not require the Palais-Smale compactness condition and provides the existence of  twin critical points.
The interaction at infinity in the electromagnetic potential implies a weak compactness property of the Poincar\'e action functional, while the interaction around the origin gives the separation property of the functional. 
}
\end{abstract}

\maketitle

\section{Introduction}

\subsection{The Lorentz force equation and the Poincar\'e action functional} 
In \cite{Fey}, The Feynman Lectures on Physics, Volume II, Chapter 26, Section 4 called 'The equations of motion in relativistic notation' 
it is written:
\begin{description}
	\item{\it "For a single charge in the fields $E$ and $B,$ the force is
		$$F = q(E + v\times B).$$
		This force is equal to the mass times the acceleration for low velocities, but the correct law for any velocity is that the force is equal to $dp/dt.$ Writing $p=m_0v/ \sqrt{1 - v^2/c^2},$ we find that the relativistically correct equation of motion is
		$$\frac{d}{dt}\left(\frac{m_0v}{\sqrt{1 - v^2/c^2}}\right) = F = q(E + v\times B)."$$ 
	}
\end{description}
Moreover, in the Introduction of his book \cite{Ebook}, Ekeland wrote:
\begin{description}
	\item{\it "The interest in periodic solution to Hamiltonian systems has been unabated for three centuries now. In an oft-quoted sentence, Poincar\'e considers them to be 'the only opening through which we can force our way into an otherwise impenetrable citadel'." }
\end{description}
In what follows the Euclidean scalar product on $\RR^3$ is denoted by
\begin{eqnarray*}
x\cdot y = x_1y_1 + x_2y_2 + x_3y_3 \quad\quad (x,y\in\RR^3),
\end{eqnarray*}
with the corresponding norm denoted by $|\cdot|.$
Let $T>0$ be a fixed period and
\begin{eqnarray*}
V:[0,T]\times \RR^3\to\RR, \quad\quad W:[0,T]\times \RR^3\to\RR^3,
\end{eqnarray*}
be two $C^1$-functions. We consider the Lorentz force equations (LFE)
\begin{eqnarray*}
\left(\frac{u'}{\sqrt{1-|u'|^2}}\right)' = E(t,u) + u'\times B(t,u), 
\end{eqnarray*}
with periodic boundary conditions 
\begin{eqnarray*}
u(0)=u(T), \quad u'(0)=u'(T),
\end{eqnarray*}
where for all $(t,\xi)\in [0,T]\times \RR^3$ one has that
\begin{eqnarray*}
E(t,\xi)=-\nabla_{\xi}V(t,\xi)-\frac{\partial W}{\partial t}(t,\xi), 
\end{eqnarray*}
and
\begin{eqnarray*}
B(t,\xi)=\mbox{curl} _{\xi} W(t,\xi).
\end{eqnarray*}
It is well known that $E+u'\times B$ is the so called {\it Lorentz force} with the {\it electric field} $E$ and the {\it magnetic field} $B.$ Moreover, $V$ is called the {\it electric scalar potential} and $W$ is the {\it magnetic vector potential.}
By a {\it $T$-periodic  solution} we understand as usual a $C^2$-function $u:[0,T]\to\RR^3$ such that $|u'(t)|<1$ for all $t\in [0,T]$ which verifies the LFE and the periodic boundary conditions on $[0,T].$


Again, in \cite{Fey}, The Feynman Lectures on Physics, Volume II, Chapter 19, Section 1 called 'The Principle of Least Action' 
it is written:
\begin{description}
	\item {\it "The function that is integrated over time to get the action $\mathbb{S}$
		is called the Lagrangian, $\mathcal L,$
		which is a function only of the velocities and positions of particles. So the principle of least action is also written
		$$\mathbb{S} = \int_{t_1}^{t_2}\mathcal L(x_i,v_i)dt,$$
		where by $x_i$ and $v_i$ are meant all the components of the positions and velocities. So if you hear someone talking about the 'Lagrangian,' you know they are talking about the function that is used to find $\mathbb{S}.$ For relativistic motion in an electromagnetic field	 
		$$\mathcal L = -m_0c^2\sqrt{1 - v^2/c^2} - q(\phi - v\cdot A)."$$ }
\end{description}
In modern terms of Nonlinear Functional Analysis the action $\mathbb{S}$ is described as follows. Consider the function spaces
\begin{eqnarray*}
W^{1,\infty}_T = \left\{u\in W^{1,\infty}(0,T;\RR^3): u(0) = u(T)\right\},
\end{eqnarray*}
and
\begin{eqnarray*}
C_T = \left\{u\in C([0,T];\RR^3): u(0) = u(T)\right\}.
\end{eqnarray*}
The norm on $L^{\infty}(0,T;\RR^3)$
is denoted by $||\cdot||_{\infty}$ and it is well-known that $(C_T, ||\cdot||_{\infty})$ and $(W^{1,\infty}_T,||\cdot||_{1,\infty})$ are Banach spaces, where
\begin{eqnarray*}
||u||_{1,\infty} = ||u||_{\infty} + ||u'||_{\infty} \quad\quad (u\in W^{1,\infty}_T).
\end{eqnarray*}
Consider the functionals 
\begin{eqnarray*}
\phi: W^{1,\infty}_T\to\RR, \quad\quad \psi: W^{1,\infty}_T\to (-\infty,+\infty]
\end{eqnarray*}
given for all $u\in W^{1,\infty}_T$ by
\begin{eqnarray*}
\phi (u) = \int_0^T\left[W(t,u)\cdot u' - V(t,u)\right]dt,
\end{eqnarray*}
and 
$$
\psi (u) = \left\{
\begin{array}{ll}
{\displaystyle \int_0^T\left[1 - \sqrt{1 - |u'|^2}\right]dt}, & \mbox{ if }u\in D(\psi),\\
\\
 +\infty,& \mbox{ otherwise,} 
\end{array}
\right. 
$$
where
\begin{eqnarray*}
D(\psi) = \left\{u\in W^{1,\infty}_T: ||u'||_{\infty}\leq 1\right\}.
\end{eqnarray*}
{\it The Poincar\'e action functional}  (see Damour \cite{Dam} for a nice introduction of Poncar\'e's contribution to LFE) associated to the Lorentz force equation with $T$-periodic boundary conditions is then given by
$$\mathcal F = \phi + \psi.$$
The key ingredient connecting the Lorentz force equation with the Poincar\'e action functional is due to Bereanu, Jebelean \& Mawhin \cite{BJM} and Arcoya, Bereanu \& Torres \cite{ABT} and reads as follows. 

\begin{theorem}\label{LFEPoincare}
	Assume that the electromagnetic potential $(V, W)$ is of class $C^1$ and $u$ is a {\bf critical point}  of the Poincar\'e action functional $\mathcal F$, that is,
	$u\in D(\psi)$ and the following variational inequality holds true,
	\begin{eqnarray*}
		\langle \phi'(u), v - u\rangle + \psi(v) - \psi(u)\geq 0, \quad\mbox{for all} \quad v\in W^{1,\infty}_T.
	\end{eqnarray*}
	Then, $u$ is a $T$-periodic solution of the Lorentz force equation. Moreover, if $W(0,\cdot) = W(T,\cdot)$ then any $T$-periodic solution of the Lorentz force equation is a critical point of the Poincar\'e action functional $\mathcal F.$
\end{theorem}

\subsection{Related results: existence of periodic solutions}
It is well known that one of the most powerful tools of Nonlinear Functional Analysis is the so called  the Direct Method of the Calculus of Variations having its roots in Hilbert's 20th problem \cite{Hilbert}: 
\begin{description}
	\item {\it "It is my conviction that it will be possible to prove these existence theorems by means of a general principle whose nature is indicated by Dirichlet's principle."}
\end{description}
After the efforts of many mathematicians among whom Lebesgue, Tonelli, Golomb, Mazur, Schauder, the above Hilbert's conviction became an elegant theorem which can be stated as follows: {\it in a reflexive Banach space any weakly lower semicontinuous functional having a bounded minimizing sequence possesses a global minimum}. The corresponding result for the Poincar\'e action functional $\mathcal F$ is given in \cite[Theorem 7]{ABT}: if there exists $\rho>0$ such that 
\begin{eqnarray*}
\inf_{\left\{u\in W^{1,\infty}_T: |\bar{u}|\leq\rho\right\}}\mathcal F = \inf_{W^{1,\infty}_T}\mathcal F,
\end{eqnarray*}
then $\mathcal F$ is bounded from below on $W^{1,\infty}_T$ and has a global minimum. Here, $\bar{u}$ denotes the mean value of $u$ on $[0,T].$ The main application of this result is given in \cite[Theorem 8]{ABT} and can be stated as follows. Let us assume that one has the following interaction at infinity between the electric potential $V$ and the magnetic potential $W:$
\begin{description}
	\item [$(VW_\infty^1)$] there exist $\mu>\nu\geq1$ and $\theta,\rho>0$ such that $V(t,\xi)\leq -\theta |\xi|^{\mu}$ and $|W(t,\xi)|\leq \theta|\xi|^{\nu}$ for all $(t,\xi)\in [0,T]\times\RR^3$ with $|\xi|\geq \rho.$
\end{description}
Then, the Lorentz force equation has at least one $T$-periodic solution which is a minimizer of the Poincar\'e action functional $\mathcal F.$ 

Moreover, it is well known that one of the most powerful minimax method used to deal with unbounded functionals (for which the direct method cannot be applied) is Rabinowitz's Saddle Point Theorem. In \cite[Corollary 2]{ABT} it is proved that this minimax method holds true also for the Poincar\'e action functional. The main application of this result is given in \cite[Theorem 9]{ABT},  and can be stated as follows. Let us assume that one has the following interaction at infinity between the electric potential $V$ and the magnetic potential $W:$
\begin{description}
	\item [$(VW_\infty^2)$] there exist $\mu>\nu\geq 1$ and $\theta,\rho>0$ such that $V(t,\xi)\geq \theta |\xi|^{\mu}$ and $|W(t,\xi)|\leq \theta|\xi|^{\nu}$ for all $(t,\xi)\in [0,T]\times\RR^3$ with $|\xi|\geq \rho.$
\end{description}
For $u\in W^{1,\infty}_T$ we denote by $\bar u$ the mean value of $u$, i.e., 
$$
\bar u=\frac{1}{T}\int_0^T u(t)dt,
$$
and we consider the decomposition
$$
W^{1,\infty}_T = \RR^3\oplus {\widetilde W^{1,\infty}_T},
$$
where
$$
{\widetilde W^{1,\infty}_T} = \{u\in W^{1,\infty}_T: \bar{u}=0\}.
$$
Thus, for any $u\in W^{1,\infty}_T$ there exists an unique element $\tilde u\in \widetilde W^{1,\infty}_T$ such that 
$$
u=\bar u+\tilde u.
$$
If $(VW_\infty^2)$ holds, then  there exists $\alpha>0$ such that 
$$
\sup_{\partial B_{\RR^3}[0,\alpha]} \mathcal F < \inf _{{\widetilde W^{1,\infty}_T}}\mathcal F,
$$
and the Lorentz force equation has at least one $T$-periodic solution $u$ which is a critical point of the Poincar\'e action functional $\mathcal F$ such that 
$$
\mathcal F(u) = \inf_{\gamma\in \Gamma} \sup_{B_{\RR^3}[0,\alpha]} \mathcal (F\circ\gamma),
$$ 
where $\Gamma = \{\gamma\in C( B_{\RR^3}[0,\alpha], W^{1,\infty}_T): \gamma = id \ \mbox{on} \ \partial B_{\RR^3}[0,\alpha]\}.$

\subsection{Main result: two periodic solutions via double interaction}

Let us introduce our interaction at infinity in the form of 
\begin{description}
	\item [$({\cal A}_{\infty})$] the electric and magnetic potentials are such that
	\begin{eqnarray*}
		\left|\int_0^T V(t,u)dt\right| - \int_0^T|W(t,u)|dt\to +\infty, \quad\mbox{as} \ u\in D(\psi), \ |\bar{u}|\to\infty.
	\end{eqnarray*}
\end{description}
We will also use two interactions around the origin, namely,
	\begin{description}
	\item [$({\cal A}_{0}^1)$] there exist $R>r>\frac{T}{4},$ $R>1+\frac{T}{2}$ 
	such that
	$$
	T\max_{[0,T]\times \left(\Gamma_{T}\cup \Gamma_R\right)} \left(V+|W|\right)< \min_{|\xi|=r}\int_0^T V(t,\xi)dt,
	$$
\end{description}
and
\begin{description}
	\item [$({\cal A}_{0}^2)$] there exist  $R>1+\frac{T}{2}$ and  $\xi_0\in\RR^3$ with $|\xi_0|>R$ 
	such that
	$$T\max_{[0,T]\times  \Gamma_R} \left(V+|W|\right)< \min_{\xi\in\{0,\xi_0\}}\int_0^T V(t,\xi)dt.$$
\end{description}
	where 
	$$
	\Gamma_{T}=B_{\RR^3}\left[0,\frac{T}{4}\right]
	$$
	and
	$$
	\Gamma_R =  \left\{\xi\in\RR^3: R-1-\frac{T}{2}\leq |\xi|\leq R\right\}.
	$$
Our main result reads as follows.
\begin{theorem}[Double interaction principle for the LFE]\label{MainRes}
Suppose that the electromagnetic potential $(V, W)$ is of class $C^1$ and  satisfies $({\cal A}_{\infty})$ together with either $({\cal A}_{0}^1)$ or $({\cal A}_{0}^2).$ Then the Lorentz force equation has at least two $T$-periodic solutions.
\end{theorem}
We have the following consequences.
\begin{corollary}\label{cor1}
	Suppose that the electromagnetic potential $(V, W)$ is of class $C^1,$ satisfies $(VW_\infty^1)$ and	
	\begin{description}
		\item [$({\cal A}_{0}^3)$] there exist $r>\frac{T}{4}$ 
		such that
		$$T\max_{[0,T]\times \Gamma_{T}} (V+|W|)< \min_{|\xi|=r}\int_0^T V(t,\xi)dt.$$
	\end{description}
	Then, the Lorentz force equation has at least two $T$-periodic solutions.
\end{corollary}

\begin{corollary}\label{cor2}
	Suppose that the electromagnetic potential $(V, W)$ is of class $C^1,$ satisfies $(VW_\infty^2)$ and $({\cal A}_{0}^1).$	
Then, the Lorentz force equation has at least two $T$-periodic solutions.
\end{corollary}

\begin{corollary}\label{cor3}
	If the electromagnetic potential $(V, W)$ is of class $C^1$ and satisfies $(VW_\infty^1)$ or $(VW_\infty^2)$  together with  $({\cal A}_{0}^2)$, then the Lorentz force equation has at least two $T$-periodic solutions.
\end{corollary}

Finally, we point out the fact that assumption $({\cal A}_{\infty})$ is a reminiscence of the well known Ahmad-Lazer-Paul condition (see \cite[Section 4]{BR}, \cite[p. 85]{MW}, \cite[p. 25]{Rabbook}). More precisely we have the following result.

\begin{corollary}\label{cor4}
		Assume the electromagnetic potential $(V, W)$ is of class $C^1$ and satisfies the following Ahmad-Lazer-Paul type condition,
		\begin{description}
			\item [$(ALP)$] $\nabla_{\xi}V$ and $W$ are bounded and 
			$$\int_0^TV(t,\xi)dt\to +\infty \  \ \mbox{as} \  \ |\xi|\to\infty,$$
			or
			$$\int_0^TV(t,\xi)dt\to -\infty \  \ \mbox{as} \  \ |\xi|\to\infty.$$
		\end{description}
		If, in addition, $(V, W)$ satisfies $({\cal A}_{0}^1)$ or $({\cal A}_{0}^2),$ then  the Lorentz force equation has at least two $T$-periodic solutions.
\end{corollary}

\subsection{Strategy of the proof}

For the proof of our main results we introduce a new minimax method for the Poincar\'e action functional. In the Introduction of \cite{Rabbook}, Rabinowitz explains the main idea of the minimax methods in general as follows:
\begin{description}
	\item {\it "What are minimax methods? These are methods that characterize a critical value $c$ of a functional $I$ as a minimax over a suitable class of sets $\mathcal S:$
		$$c = \inf_{A\in\mathcal S} \max_{u\in A} I(u) \quad \quad (*)$$
		There is no recipe for choosing $\mathcal S.$ In any given situation the choice must reflect some qualitative change in the topological nature of the level sets of $I,$ i.e. $I^{-1}(s)$ for $s$ near $c.$ Thus obtaining and characterizing a critical value $c$ as in $(*)$ is something of an ad hoc process."
	}
\end{description}
In the same line Nirenberg explains in \cite{Nir}, Section II.2:
\begin{description}
	\item {\it "We shall present several arguments for finding stationary points of a
		functional $f$ which is unbounded from above and below - leading us to look
		for saddle points. These are found by some min max argument. One considers
		a certain class $\mathcal C$ of sets $\Sigma$ and forms 
		$$\max_{u\in\Sigma}f(u)$$
		for any $\Sigma$ in $\mathcal C$ and then 
		$$c = \inf_{\Sigma\in\mathcal C} \max_{u\in\Sigma}f(u),$$
		and tries to show that the number $c$ is a stationary value of $f$ i.e., is a value of $f$ at some stationary point." 
	}
\end{description}

Concerning the Poincar\'e action functional, one has that $\mathcal F$ splits as $\mathcal F = \phi + \psi,$ with $\phi$ being $C^1$ and $\psi$ convex lower semi-continuous. To the best of our knowledge the first paper dealing with minimax methods (more precisely, the Lusternik-Schnirelmann method for $S^1$-invariant functionals) for convex perturbations of smooth functionals is the Ekeland \& Lasry article \cite{EL}. Hence, following \cite[p. 307]{EL}, 
\begin{eqnarray*}
	\partial \mathcal F(u) =  \phi'(u) + \partial\psi(u), \quad\mbox{for all} \ u\in D(\psi),
\end{eqnarray*}
which implies that for $u\in D(\psi)$ it holds that $0\in 	\partial \mathcal F(u)$ iff $-\phi'(u)\in \partial\psi(u)$ iff $u$ is a critical point of $\mathcal F$ (see Theorem \ref{LFEPoincare}). 

Next, we introduce our main tool. We define $\mathcal H$ to be the class of mappings $\gamma\in C([0,1]\times W^{1,\infty}_T, W^{1,\infty}_T)$ with the following properties:
\begin{description}
	\item [$(a)$] $\gamma(0,\cdot)=id;$ 
	
	\item [$(b)$] for each $t\in [0,1)$ one has that $\gamma(t,\cdot)$ is a homeomorphism of $W^{1,\infty}_T$ onto itself and $\gamma(\cdot,\cdot)^{-1}$ is continuous on $[0,1)\times W^{1,\infty}_T;$
	
	\item [$(c)$] there exists $u_0\in W^{1,\infty}_T$ with $\gamma(1,u)=u_0$ for all $u\in W^{1,\infty}_T$ and $\gamma(t,u)\to u_0$ as $t\to 1$ uniformly on bounded subsets of $W^{1,\infty}_T.$  
\end{description}
Let $A,B$ be two nonempty subsets of $W^{1,\infty}_T$. Following Schechter \cite{Schechter-Book99}, we say $A$ {\it links} $B$ with respect to ${\cal H}$  if 
$A\cap B = \varnothing$ and for every $\gamma\in\mathcal H$ there exist $t_A\in (0,1]$ such that 
$$
\gamma(t_A,A)\cap B\neq\varnothing.
$$
If $A$ links $B$ and $B$ links $A$ with respect to $\cal H$, then we say they form a double linking pair.
For applications  of the above linking in the $C^1$ context see for example \cite{Sajm}, \cite{Scpam}, \cite{Sduke}, \cite{Schechter-Zou-JFA03}.

\begin{theorem}[Minimax theorem for the Poincar\'e action functional]\label{MinMaxThm}
Suppose the electromagnetic potential $(V,W)$ is of class $C^1$ and the following condition
		\begin{description}
			\item [(C)] $ |\phi(u)|\to\infty, \quad \mbox{as} \ u\in D(\psi), \ |\bar{u}|\to\infty$,
		\end{description}	
holds together with the separation property			
		\begin{description}	
			\item [(S)] there exist two nonempty, closed and bounded subsets $A$ and $B$ of $W^{1,\infty}_T$ forming a double linking pair such that $A\subset D(\psi)$ and
			$$a=\sup_A \mathcal F \leq \inf_B \mathcal F=b.$$
		\end{description}
Then $\mathcal F$ possesses two critical values c, d such that $c\geq b$ and $d\leq a$ given by
		$$c= \inf_{\gamma\in\mathcal H}\sup_{[0,1]\times A} (\mathcal F\circ\gamma), \quad\quad d=\sup_{\gamma\in\mathcal H}\inf_{[0,1]\times B} (\mathcal F\circ\gamma).$$
		Moreover, if $c=b$ then $\mathcal F$ possesses a critical point in $\overline{B}^{C_T}$ at the level $c$, whereas $d=a$ ensures that $\mathcal F$ possesses a critical point in $\overline{A}^{C_T}$ at the level $d.$	In particular, the two critical points are distinct provided $c\neq d$ or $\overline{A}^{C_T}\cap \overline{B}^{C_T}\neq \varnothing.$
\end{theorem}
We point out the fact that condition $(C)$ is the only compactness input we need throughout the paper. As it will be seen in the subsequent sections, this condition is sufficient for any Palais-Smale sequence of the functional ${\cal F}$ to be bounded and possess a  subsequence that converges in $C_T$ to a critical point of $\F$.

In order to prove our main result, namely Theorem \ref{MainRes}, we proceed as follows: if assumption $({\cal A}_0^1)$ holds true, then we take 
$$
A = \left\{\xi\in\RR^3: |\xi| = r\right\}
$$
and
$$
B = \left\{u\in{\widetilde W^{1,\infty}_T}: ||u||_{1,\infty}\leq R \right\}\cup \left\{u+se: u\in{\widetilde W^{1,\infty}_T}, \ s\geq 0, \ ||u+se||_{1,\infty} = R\right\},
$$
where $e\in\RR^3$ with $|e|=1$ is fixed. Moreover, when assumption $({\cal A}_0^2)$ holds true we take
$$A=\{0, \xi_0\}, \quad B=\partial B_{W^{1,\infty}_T}[0,R],
$$
then we show that $(S)$ is satisfied with $A, B$ defined above. On the other hand we will prove that condition $(C)$  is implied by the interaction at infinity $({\cal A}_{\infty}).$

\section{Qualitative properties of the Poincar\'e action functional}

In what follows we use the notation
\begin{eqnarray*}
	(X, ||\cdot||_X) = ( W^{1,\infty}_T, ||\cdot||_{1,\infty}),\ (Y, ||\cdot||_Y) = (C_T, ||\cdot||_{\infty}) \mbox{ and } \tilde X=\widetilde W^{1,\infty}_T.
\end{eqnarray*}

We start with an embedding result which plays a key role in the proof of the minimax theorem presented in the previous section.

\begin{lemma}\label{CompEmbed}
The embedding $X\hookrightarrow Y$ is compact. Moreover, if $\{u_n\}\subset  X$ is a bounded sequence in $X$ such that $u_n\to u$ in $Y,$ then $u\in  X$ and 
	$$
	\liminf_{n\to \infty} \|u_n\|_X\geq \|u\|_X.
	$$
\end{lemma}
\begin{proof}
A bounded subset of $X$ is bounded and equi-Lipschitz in $Y$, hence relatively compact by the Ascoli-Arzel\`a theorem, so the embedding $X\hookrightarrow Y$ is indeed compact.  The fact that $u\in X$ and that $u_n'\rightharpoonup^*u'$ in $\sigma(L^\infty,L^1)$ is proved in  \cite[Lemma 12]{Maw11} (see also \cite{BJM} ). Only the lower semicontinuity of the norm has to be added, and it follows from the weak$^*$ lower semicontinuity of $\|\cdot\|_{L^\infty}$ together with $\|u\|_Y=\lim_{n\to\infty}\|u_n\|_Y$:
\begin{eqnarray*}
\|u\|_X &=& \|u\|_Y+\|u'\|_Y\le \liminf_{n\to\infty}\|u_n\|_Y+\liminf_{n\to\infty}\|u_n'\|_Y\\
	&\le& \liminf_{n\to\infty}\left(\|u_n\|_Y+\|u_n'\|_Y\right)=\liminf_{n\to\infty}\|u_n\|_X.
\end{eqnarray*}
\end{proof}

Concerning the nonsmooth part $\psi$ of the Poincar\'e action $\mathcal F$ we need the following basic estimates. 
\begin{lemma}\label{basicineq}
	Let $u\in D(\psi)$ be fixed. The following assertions hold true:
	\begin{description}
		\item [(i)] $||\tilde{u}||_{Y}\leq\frac{T}{4}||u'||_{Y}\leq \frac{T}{4}.$
		In particular, $u([0,T])\subset B_{\RR^3}\left [\bar{u}, \frac{T}{4}\right]$;
		\item [(ii)] $||u||_{Y}\geq ||u||_{X}-1$;
		\item [(iii)] $\psi(u)\geq \frac{2}{T}\|\tilde u\|_Y^2$;
		\item [(iv)] if $\bar{u}=0,$ then $||u||_{Y}\leq \frac{T}{4+T}||u||_{X}$;
		\item [(v)] $|u(t)|\geq ||u||_{X}-1-\frac{T}{2}$, for all $t\in [0,T]$. 
	\end{description}
\end{lemma}

\begin{proof}
\begin{description}

\item (i)  For $0\le s\le t\le T$ one has 
$$
|u(t)-u(s)|=\left|\int_s^tu'(\tau)d\tau\right|\le\|u'\|_Y(t-s).
$$ 
Using the above estimate and the fact that  $u(T)=u(0)$ we get
\begin{eqnarray*}
|u(t)-u(s)| & = & \left|\left(u(t)-u(T)\right)+\left(u(0)-u(s)\right)\right|\\
	&\le &\|u'\|_Y(T-t)+\|u'\|_Y s\\
	& = &\|u'\|_Y\left(T-(t-s)\right). 
\end{eqnarray*}
Consequently,
\begin{equation}\label{CircLip}
|u(t)-u(s)|\le\|u'\|_Y\min\left\{|t-s|,T-|t-s|\right\},\quad \forall s,t\in[0,T].
\end{equation}
We show next that 
\begin{equation}\label{MinInt}
\int_0^T\min\left\{|t-s|,T-|t-s|\right\}ds=\frac{T^2}{4}, \quad \forall t\in [0,T].
\end{equation}
Let us first fix $t\in [0, T/2]$. Then 
$$
\min\left\{|t-s|,T-|t-s|\right\}=\left\{
\begin{array}{ll}
t-s, & \mbox{ if }0\leq s\leq t,\\
s-t, & \mbox{ if } t<s\leq t+\frac{T}{2},\\
T-s+t, & \mbox{ if } t+\frac{T}{2}<s\leq T,
\end{array}
\right.
$$
hence
\begin{eqnarray*}
\int_0^T\min\left\{|t-s|,T-|t-s|\right\}ds &=& \underbrace{\int_0^t (t-s)ds}_{=I_1}+\underbrace{\int_t^{t+T/2} (s-t)ds}_{=I_2}\\
	&&+\underbrace{\int_{t+T/2}^T (T-s+t)ds}_{=I_3}.
\end{eqnarray*}
Direct computations yield 
$$
I_1+I_2=\frac{t^2}{2}+\frac{T^2}{8},
$$
whereas for $I_3$, performing the substitution $\tau=T-s+t$, we immediately get 
$$
I_3=\frac{T^2}{8}-\frac{t^2}{2},
$$
hence \eqref{MinInt} holds.

If $t\in [T/2,T]$, then $T-t\in [0,T/2]$ and reasoning as above we get the desired estimate.

\smallskip

Since $\int_0^T\tilde u(s)\,ds=0$, one has $\tilde u(t)=\frac{1}{T}\int_0^T\left(u(t)-u(s)\right)ds$, whence, 
by \eqref{CircLip} and \eqref{MinInt} one has
\begin{eqnarray*}
|\tilde u(t)| & \le & \frac{1}{T}\int_0^T|u(t)-u(s)|\,ds\le\frac{\|u'\|_Y}{T}\int_0^T\min\left\{|t-s|, T-|t-s|\right\}ds\\
	& = & \frac{\|u'\|_Y}{T}\cdot\frac{T^2}{4}=\frac{T}{4}\|u'\|_Y\le\frac{T}{4},
\end{eqnarray*}
for every $t\in [0,T]$, and $|u(t)-\bar u|=|\tilde u(t)|$ gives the inclusion $u([0,T])\subset B_{\RR^3}\left[\bar u,\frac{T}{4}\right]$.

\medskip

\item (ii) By the definition of the norm of $X$ one has 
$$
\|u\|_Y=\|u\|_X-\|u'\|_Y\geq \|u\|_X-1,
$$ 
due to the fact that $u\in D(\psi)$.

\medskip

\item (iii)  Fix $s,t\in[0,T]$ with $s\le t$. The two arcs of the closed curve $u$ joining $u(s)$ to $u(t)$ have lengths $\int_s^t|u'(\tau)|d\tau$ and $\ell-\int_s^t|u'(\tau)|d\tau$, where $\ell$ is the length of the entire curve. Then $|u(t)-u(s)|$ is at most the smaller of the two, hence at most $\ell/2$. Averaging in $s$,
\begin{equation}\label{CLength}
|\tilde u(t)|=\left|\frac{1}{T}\int_0^T\left(u(t)-u(s)\right)ds\right|\le\frac{1}{T}\int_0^T\frac{\ell}{2}\,ds=\frac{\ell}{2}.
\end{equation}
We show next that 
\begin{equation}\label{PsiEst}
\psi(u)\ge\frac{1}{2}\int_0^T|u'(t)|^2dt.
\end{equation} 
For $\sigma\in[0,1]$ one has $\sqrt{1-\sigma^2}\le1-\frac{\sigma^2}{2}$, since squaring the two nonnegative sides gives $1-\sigma^2\le1-\sigma^2+\frac{\sigma^4}{4}$. Therefore $1-\sqrt{1-\sigma^2}\ge\sigma^2/2$ and, applying this with $\sigma=|u'(t)|\le1$ and integrating,
$$
\psi(u)=\int_0^T\left[1-\sqrt{1-|u'(t)|^2}\right]dt\ \ge\ \frac{1}{2}\int_0^T|u'(t)|^2dt .
$$
By the H\"{o}lder inequality, 
$$
\ell=\int_0^T|u'(t)|\cdot1\,dt\le\left(\int_0^T|u'(t)|^2dt\right)^{1/2}\sqrt T,
$$
that is $\int_0^T|u'(t)|^2dt\ge\ell^2/T$. Relations \eqref{CLength} and \eqref{PsiEst} now give
$$
\psi(u) \ge \frac{\ell^2}{2T}\ge \frac{(2\|\tilde u\|_Y)^2}{2T} =\frac{2}{T}\|\tilde u\|_Y^2 .
$$

\medskip

\item (iv) If $\bar u=0$ then $\tilde u=u$, so $\|u\|_Y=\|\tilde u\|_Y$. Now (i) gives
$$
\|u\|_Y \le \frac{T}{4}\|u'\|_Y= \frac{T}{4}\left(\|u\|_X-\|u\|_Y\right).
$$
Bringing $\|u\|_Y$ to the left-hand side, we have $\left(1+\frac T4\right)\|u\|_Y\le\frac{T}{4}\|u\|_X$, i.e.,
$$
\|u\|_Y \le \frac{\frac{T}{4}}{1+\frac{T}{4}}\|u\|_X =\frac{T}{4+T}\|u\|_X.
$$

\medskip
\item (v) For every $t\in [0,T]$, one has 
$$
|u(t)|\le|\bar u|+|\tilde u(t)|\le|\bar u|+\|\tilde u\|_Y.
$$
Taking the supremum as $t\in [0,T]$ and using (i), we have
$$
|\bar u| \ge \|u\|_Y-\|\tilde u\|_Y\ge\|u\|_Y-\frac{T}{4}.
$$
Hence, again by (i) and then by (ii),
$$
|u(t)| \ge |\bar u|-|\tilde u(t)| \ge\left(\|u\|_Y-\frac{T}{4}\right)-\frac{T}{4}=\|u\|_Y-\frac{T}{2} \ge \|u\|_X-1-\frac T2
$$
for every $t\in[0,T]$ (the inequality being automatically fulfilled when $\|u\|_X\le1+\frac T2$).
\end{description}
\end{proof}

The main functional analytic properties of the Poincar\'e action $\mathcal F$ are summarized in the following result.

\begin{lemma}\label{prop1234}
Assume that the electromagnetic potential $(V, W)$ is of class $C^1.$ Then the Poincar\'{e} action functional $\cal F$ enjoys the following properties:
\begin{description}
	\item $(P_1)$ The functional $\phi$ is $C^1$ on $X$ and for all $u, v\in  X$ one has that 
	\begin{eqnarray*}
	\langle \phi'(u), v \rangle = \int_0^T [\mathcal E(t,u,u')-\nabla_uV(t,u)]\cdot v dt + \int_0^T W(t,u)\cdot v' dt,
	\end{eqnarray*}
	where, for all $(t,\xi,\zeta)\in [0,T]\times\RR^3\times\RR^3$,
	\begin{eqnarray*}
	\mathcal E(t,\xi,\zeta)=\left(\zeta\cdot D_{\xi_1}W(t,\xi), \zeta\cdot D_{\xi_2}W(t,\xi), \zeta\cdot D_{\xi_3}W(t,\xi)\right).
	\end{eqnarray*}
	
	\item $(P_2)$ If $\{u_n\}\subset  X$ is a bounded sequence in $X$ such that $u_n\to u$ in $Y,$ then
	\begin{eqnarray*}
&	\phi(u_n)\to\phi(u), \\	&||\phi'(u_n) - \phi'(u)||_{X^*}\to 0, \\
&	\langle \phi'(u_n), u_n\rangle \to 	\langle \phi'(u), u\rangle.
	\end{eqnarray*}
	
	\item $(P_3)$ $\phi$ and $\phi'$ map bounded subsets of $X$ to bounded sets, i.e.,  $\sup_{\|u\|\le r}|\phi(u)|<+\infty$ 	and $\sup_{\|u\|\le r}\|\phi'(u)\|_{X^*}<+\infty$ for every $r>0$.
	
	\item $(P_4)$ For every $r>0$, $\phi'$ is uniformly continuous from $( B_X[0,r],\|\cdot\|_Y)$ to $(X^*,\|\cdot\|_{X^*})	$, 	i.e., for each $\varepsilon>0$ there is $\delta>0$ such that 
		$$
		\left.
		\begin{array}{l}
		u_1,u_2\in B_X[0,r]\\
		\|u_1-u_2\|_Y\le\delta
		\end{array}
		\right]\Longrightarrow \|\phi'(u_1)-\phi'(u_2)\|_{X^*}\le\varepsilon.
		$$

	\item $(P_5)$ The functional $\psi$ is proper, convex and satisfies 
	\begin{eqnarray*}
	0\leq\psi(u)\leq T, \quad\mbox{for all } u\in D(\psi).
	\end{eqnarray*}
	Moreover, if $\{u_n\}\subset  X$ is a bounded sequence in $X$ such that $u_n\to u$ in $Y,$ then
	\begin{eqnarray*}
	\psi(u)\leq\liminf _{n\to\infty}\psi(u_n).
	\end{eqnarray*}
	
	\item $(P_6)$ Let $A, B$ be two bounded subsets in $X$ such that $A\subset D(\psi).$ Then, 	there exists 	  			$\tilde{\gamma}\in\mathcal H$ such that
	\begin{eqnarray*}
		\sup_{[0,1]\times A} (\mathcal F\circ \tilde{\gamma})<+\infty \quad \mbox{and} \quad \inf_{[0,1]\times B} (\mathcal F\circ \tilde{\gamma})>-\infty.
	\end{eqnarray*}
\end{description}
\end{lemma}

\begin{proof} 

\begin{description}
\item $(P_1)$ See \cite[Lemma 1]{ABT}.

\item $(P_2)$ The first and the third convergence were established in \cite[p. 1705]{ABT}, hence we only need to prove the second convergence. 

Let $\{u_n\}\subset X$ be a bounded sequence such that $u_n\to u$ in $Y$. Lemma \ref{CompEmbed} ensures $u\in X$ and there exists $r>0$ such that 
$$
\|u\|_X\le r \mbox{ and } \|u_n'-u\|_Y\le 2r, \quad \forall n \ge 1.
$$
Let $\mathcal M(t,\xi)$ be the $3\times3$ matrix whose columns are $D_{\xi_1}W(t,\xi),D_{\xi_2}W(t,\xi),D_{\xi_3}W(t,\xi)$. Then 
\begin{equation}\label{Mmatrix}
{\cal E}(t,\xi,\zeta)=\mathcal M(t,\xi)^{\top}\zeta,\quad {\cal E}(t,\xi,\zeta)\cdot\eta=\zeta\cdot\mathcal M(t,\xi)\eta,\quad\forall \zeta,\eta\in\mathbb{R}^3,
\end{equation}
and put $\mathcal M(t)=\mathcal M(t,u(t))\in C([0,T];\mathbb R^{3\times3})$ and $c_K=\|\mathcal M\|_\infty$. Since $\nabla_uV$, $W$ and $\mathcal M$ are uniformly continuous on the compact set $[0,T]\times B_{\RR^3}[0,r]$, we infer that the sequence 
\begin{eqnarray*}
\omega_n &=&  \left\|\nabla_uV(\cdot,u_n)-\nabla_uV(\cdot,u)\right\|_\infty+\left\|W(\cdot,u_n)-W(\cdot,u)\right\|_\infty\\
	&& +\left\|\mathcal M(\cdot,u_n)-\mathcal M(\cdot,u)\right\|_\infty
\end{eqnarray*}
converges to $0$ as $n\to\infty.$

Let $v\in X$ with $\|v\|_X\le1$, so that $\|v\|_Y\le1$ and $\|v'\|_Y\le1$. By \eqref{Mmatrix},
$$
{\cal E}(t,u_n,u_n')-{\cal E}(t,u,u')=\left[\mathcal M(t,u_n)-\mathcal M(t)\right]^{\top}u_n'+\mathcal M(t)^{\top}(u_n'-u'),
$$
and $\langle\phi'(u_n)-\phi'(u),v\rangle$ splits into the four integrals 
$$
I_1=\int_0^T\left[\mathcal M(t,u_n)-\mathcal M(t)\right]^{\top}u_n'\cdot v\,dt,\quad
I_2=\int_0^T(u_n'-u')\cdot\mathcal M(t) vdt,
$$
$$
I_3=-\int_0^T\left[\nabla_uV(t,u_n)-\nabla_uV(t,u)\right]\cdot vdt,\quad
I_4=\int_0^T\left[W(t,u_n)-W(t,u)\right]\cdot v'dt.
$$
It is readily seen that  $I_1\le rT\omega_n$ and $I_3,I_4\le T\omega_n$ respectively, uniformly in $v$, therefore tend to $0$.

It remains to show that $I_2\to 0$ as $n\to \infty$. Let $\varepsilon\in(0,1]$ and let $\mathcal N:[0,T]\to\mathbb R^{3\times3}$ be of class $C^1$ with $\|\mathcal M-\mathcal N\|_\infty\le\varepsilon$, which is possible because the $C^1$ functions are dense in $C([0,T];\mathbb R^{3\times3})$. Then $\|\mathcal N\|_\infty\le c_K+1$ and, $v$ being Lipschitz, so is $\mathcal N v$, with
$$
\|\mathcal N v\|_\infty\le c_K+1,\quad
\left\|(\mathcal Nv)'\right\|_\infty=\left\|\mathcal N'v+\mathcal N v'\right\|_\infty\le\|\mathcal N'\|_\infty+c_K+1=:c_\varepsilon ,
$$
by $\|v\|_\infty\le1$ and $\|v'\|_\infty\le1$. Thus
\begin{eqnarray*}
I_2 &=& \int_0^T(u_n'-u')\cdot(\mathcal M-\mathcal N+\mathcal N) vdt\\
	&=&	\int_0^T(u_n'-u')\cdot(\mathcal M-\mathcal N)vdt +\int_0^T (u_n'-u')\cdot\mathcal Nvdt\\
	&=&\int_0^T(u_n'-u')\cdot(\mathcal M-\mathcal N)vdt+\left[(u_n-u)\cdot\mathcal N v\right]\big|_0^T
-\int_0^T(u_n-u)\cdot(\mathcal Nv)'\,dt\\
	&\le& 2rT\varepsilon+ 2c_\varepsilon\|u_n-u\|_Y+Tc_\varepsilon\|u_n-u\|_Y\to 0,\quad \mbox{ as }\varepsilon \to 0 \mbox{ and }n\to \infty.
\end{eqnarray*}
Consequently
$$
\|\phi'(u_n)-\phi'(u)\|_{X^*}=\sup_{\|v\|\le1}\left|\langle\phi'(u_n)-\phi'(u),v\rangle\right|\to0. 
$$

\item $(P_3)$ If $\|\phi'(u_n)\|_{X^*}\to\infty$ for some $u_n\in B_X[0,r]$, then by Lemma \ref{CompEmbed} a subsequence satisfies $u_{n_k}\to w$ in $Y$ with $w\in X$. Using property $(P_2)$ we infer that $\phi'(u_{n_k})\to\phi'(w)$ in $X^*$, so $\|\phi'(u_{n_k})\|_{X^*}$ is bounded, a contradiction. The argument for $\phi$ is identical.

\item $(P_4)$ Arguing by contradiction, assume  there exist $\varepsilon>0$ and two sequences $\{u_n^1\},\ \{u_n^2\}\subset B_X[0,r]$ with $\|u_n^1-u_n^2\|_Y\to0$, but $\|\phi'(u_n^1)-\phi'(u_n^2)\|_{X^*}\ge\varepsilon$. By Lemma \ref{CompEmbed}, up to a subsequence,  $u_n^1\to w$ in $Y$ with $w\in X$. Then also $u_n^2\to w$ in $Y$. Now property $(P_2)$ ensures that $\phi'(u_n^1)\to\phi'(w)$ and $\phi'(u_n^2)\to\phi'(w)$ in $X^*$, so $\|\phi'(u_n^1)-\phi'(u_n^2)\|_{X^*}\to0$, which is clearly a contradiction.

\item $(P_5)$ See \cite[p. 149]{Maw11}.

\item $(P_6)$ Let $\rho>0$ be such that $A\cup B\subset B_{X}[0,\rho]$ and $$\theta = \max_{[0,T]\times B_{\RR^3}[0,\rho]}(|V|+|W|).$$ Consider $ \tilde{\gamma}\in C([0,1]\times X, X)$ given by 
\begin{eqnarray*}
	\tilde\gamma (t,u) = (1-t)u, \quad\mbox{for all} \ (t,u)\in [0,1]\times X.
\end{eqnarray*}
It follows that $\tilde\gamma (0,\cdot) = id$ and for every $t\in [0,1)$ one has that $\tilde\gamma (t,\cdot)$ is a homeomorphism with $\tilde\gamma (t,u)^{-1} = \frac{1}{1-t}u$ for all  $u\in X.$ Moreover, $\tilde\gamma (1,\cdot) = 0$ and it is clear that $\tilde\gamma (t,u)\to 0$ as $t\to 1$ uniformly on bounded subsets of $X.$ Hence, $\tilde\gamma \in\mathcal H.$ Next, it is readily seen that $(1-t)u\in D(\psi)$ and $|(1-t)u(t)|\leq \rho$ for all $t\in [0,1]$ and $u\in A,$ so from  $0\leq \psi\leq T$ we deduce that 
$$\sup_{[0,1]\times A} (\mathcal F\circ \tilde\gamma )\leq T+T\theta<+\infty.$$
Analogously, we have that 
$$\inf_{[0,1]\times B} (\mathcal F\circ \tilde\gamma )\geq -T\theta>-\infty,$$
and the proof is completed. 

\end{description}
\end{proof}

\section{Proof of the main results}

\subsection{The Double Interaction Principle}

We start by showing that a good interaction around the origin in the electromagnetic potential $(V, W)$ ensures a good geometry of the Poincar\'e action $\mathcal F$.

\begin{lemma}\label{a01S}
	Assume the electromagnetic potential $(V, W)$ satisfies $({\cal A}_0^1)$  and consider 
	$$
	A = \partial B_{\RR^3}(0,r)
	$$
	and
	$$
	B = \left\{u\in{\tilde X}: ||u||_{X}\leq R \right\}\cup \left\{u+se: u\in{\tilde X}, \ s\geq 0, \ ||u+se||_{X} = R\right\},
	$$
	where $e\in\RR^3$ with $|e|=1$ is fixed and $0<r<R$. Then, the Poincar\'e action $\mathcal F$ satisfies
	$$\sup_A \mathcal F < \inf_B \mathcal F<+\infty.$$
\end{lemma}
\begin{proof} Consider $\xi\in A,$ so $\xi\in\RR^3$ with $|\xi| = r.$ It follows that $$\mathcal F(\xi)=-\int_{0}^{T} V(t,\xi)dt,$$
and
$$\sup_A \mathcal F = \max_{|\xi|=r}\left(-\int_{0}^{T} V(t,\xi)dt\right).$$
Let $u\in B\cap D(\psi)$ be fixed. We claim that 
$$u([0,T])\subset \Gamma_T\cup \Gamma_R.$$
Indeed, if $\bar{u} = 0,$ then using Lemma \ref{basicineq} (i) we deduce that $||u||_{Y}\leq \frac{T}{4}$ and $u([0,T])\subset \Gamma_T.$ If $\bar{u} \neq 0,$ then $u\notin {\tilde X}$ and using that $u\in B$ it follows that $u = z+se$ with $z\in {\tilde X},$ $s>0$ and $||u||_{X} = R.$
Hence, using Lemma \ref{basicineq} (v) we deduce that for all $t\in [0,T],$
$|u(t)|\geq ||u||_{X} - 1 - \frac{T}{2} = R - 1 - \frac{T}{2}$ and $|u(t)|\leq ||u||_{Y}\leq ||u||_{X} = R.$ So, $u([0,T])\subset \Gamma_R,$ and our claim is proved. On the other hand, from  $\psi\geq 0$ we deduce that 
\begin{eqnarray*}
	\mathcal F(u) &=& \psi(u)  + \int_0^T\left[W(t,u)\cdot u' - V(t,u)\right]dt \\
	&\geq& -\int_0^T |W(t,u)||u'|dt- \int_0^T V(t,u)dt \\
	&\geq& -\int_0^T[|W(t,u)|+V(t,u)]dt\\
	&\geq& -T\max_{[0,T]\times (\Gamma_T\cup \Gamma_R)} (V+|W|). 
\end{eqnarray*}
So, taking into account that $\mathcal F(v) = +\infty$ for all $v\in X\setminus D(\psi)$ and $({\cal A}_0^1)$ we deduce that 
$$\inf_B \mathcal F\geq -T\max_{[0,T]\times (\Gamma_T\cup \Gamma_R)} (V+|W|)>-\min_{|\xi|=r}\int_0^T V(t,\xi)dt = \sup_A \mathcal F.$$
Moreover, one has that $\inf_B \mathcal F<+\infty$ because $0\in B\cap D(\psi).$
\end{proof}

\begin{lemma}\label{a02S}
		Assume $({\cal A}_0^2)$ to be  satisfied and consider 
		$$A=\{0, \xi_0\}, \quad B=\partial B_{X}[0,R].$$
		Then,
		$$\sup_A \mathcal F < \inf_B \mathcal F<+\infty.$$
\end{lemma}
\begin{proof} Again, for all  $\xi\in\RR^3$ one has that $$\mathcal F(\xi)=-\int_{0}^{T} V(t,\xi)dt,$$
so
$$\sup_A \mathcal F = - \min_{\xi\in\{0,\xi_0\}}\int_0^T V(t,\xi)dt.$$
Let $u\in B\cap D(\psi)$ be fixed. Using Lemma \ref{basicineq} (v) it follows that for all $t\in [0,T],$
$$
R=||u||_{X}\geq |u(t)|\geq ||u||_{X} - 1 - \frac{T}{2}=R-1-\frac{T}{2},
$$
and $u([0,T])\subset  \Gamma_R.$ Since $\psi\geq 0$, we deduce that 
\begin{eqnarray*}
	\mathcal F(u) &=& \psi(u)  + \int_0^T[W(t,u)\cdot u' - V(t,u)]dt \\
	&\geq& -\left|\int_0^T W(t,u)\cdot u'dt\right| - \int_0^T V(t,u)dt \\
	&\geq& -T\max_{[0,T]\times  \Gamma_R} (V+|W|). 
\end{eqnarray*}
So, taking into account that $\mathcal F(v) = +\infty$ for all $v\in X\setminus D(\psi)$ and $({\cal A}_0^2)$ we deduce that 
$$
\inf_B \mathcal F\geq -T\max_{[0,T]\times  \Gamma_R} (V+|W|) > \sup_A \mathcal F.
$$
Moreover, one has that $\inf_B \mathcal F<+\infty$ because $\partial B_{\RR^3}[0,R]\subset B\cap D(\psi)$ and 
$$
\inf_B \F\leq \mathcal F(\xi) = -\int_0^T V(t,\xi)dt\le -T\min_{[0,T]\times \partial B_{\RR^3}[0,R]}V<+\infty.
$$
\end{proof}

Now we show that a good interaction at infinity in the electromagnetic potential $(V, W)$ ensures the compactness property of the Poincar\'e action $\mathcal F.$

\begin{lemma}\label{AinfinitC}
	Suppose that $(V, W)$ satisfies $({\cal A}_{\infty}).$ Then, $\mathcal F$ satisfies condition $(C).$ 
\end{lemma}
\begin{proof} If $u\in D(\psi),$ then $||u'||_Y \leq 1$ and
\begin{eqnarray*}
	\left|\int_0^T V(t,u)dt\right| - |\phi(u)|&\leq& \left|\phi(u) + \int_0^T V(t,u)dt\right| \\
	&=& \left|\int_0^T W(t,u)\cdot u'dt\right| \\
	&\leq & \int_0^T |W(t,u)dt|,
\end{eqnarray*}
hence,
\begin{eqnarray*}
		\left|\int_0^T V(t,u)dt\right| - \int_0^T |W(t,u)dt| \leq |\phi(u)|,
\end{eqnarray*}
and $({\cal A}_{\infty})$ implies $(C).$
\end{proof}

\medskip

\noindent {\it Proof of Theorem \ref{MainRes}.} 
Assume $({\cal A}_{\infty})$ and $({\cal A}_0^1)$ hold true. From Lemma \ref{AinfinitC} it follows that condition $(C)$ is satisfied. Using Lemma \ref{a01S} it follows that $(S)$ is satisfied with $A, B$ defined in Lemma \ref{a01S}. Note that from \cite[Example 7, Section 2.6]{Schechter-Book99} the sets $A, B$ form a double linking pair. Now,  Theorem \ref{MinMaxThm} ensures that there exist $u$ and $v$ two critical points of $\mathcal F$ such that 
$$
\mathcal F(u) = d\leq \sup_A \mathcal F < \inf_B \mathcal F \leq c=\mathcal F(v),
$$
hence, $u\neq v.$ Moreover, using Theorem \ref{LFEPoincare} it follows that $u, v$ are $T$-periodic solutions of the Lorentz force equation.

If $({\cal A}_{\infty})$ and $({\cal A}_0^2)$ hold true, then the proof is the same but with Lemma \ref{a01S} replaced by Lemma \ref{a02S}. Notice that $A, B$ defined in Lemma \ref{a02S} is a double linking pair by \cite[Example 1, Section 2.6]{Schechter-Book99}. 
\qed

\subsection{The Corollaries}

\begin{lemma}\label{a01a02ainfinit}
	If the electromagnetic potential satisfies $(VW_\infty^1)$ or $(VW_\infty^2)$, then $({\cal A}_{\infty})$ holds true. 
\end{lemma}
\begin{proof} 	Assume that $(V, W)$ satisfies $(VW_\infty^2).$ Then there exists $c_1>0$ such that for all $(t,\xi)\in [0,T]\times \RR^3,$
\begin{eqnarray*}
	V(t,\xi)\geq \theta |\xi|^{\mu} + c_1, \quad |W(t,\xi)|\leq \theta |\xi|^{\nu} +c_1.
\end{eqnarray*}
Let $u\in D(\psi)$ be fixed, so $||u'||_Y \leq 1$ and $||\tilde{u}||_Y\leq \frac{T}{4}<T.$
It follows that 
\begin{eqnarray*}
	\int_0^T |W(t,u)|dt &\leq& \theta \int_0^T |u|^{\nu}dt + Tc_1 \\
	&\leq& 2^{\nu -1}\theta \int_0^T (|\tilde{u}|^{\nu} + |\bar{u}|^{\nu})dt + Tc_1 \\
		&\leq& 2^{\nu -1}\theta  (T^{\nu + 1} + T|\bar{u}|^{\nu}) + Tc_1.
\end{eqnarray*}
Next, one has the following calculation,
\begin{eqnarray*}
	T|\bar{u}|^{\mu} &\leq& \int_0^T (|u| + |\tilde{u}|)^{\mu}dt\\
	&\leq& 2^{\mu - 1} \int_0^T (|u|^{\mu} + |\tilde{u}|^{\mu})dt \\
		&\leq& 2^{\mu - 1} \int_0^T (|u|^{\mu} +T^{\mu})dt,
\end{eqnarray*}
so, there exist $c_2, c_3>0$ two constants such that 
\begin{eqnarray*}
	\int_0^T |u|^{\mu}dt\geq c_2|\bar{u}|^{\mu} - c_3.
\end{eqnarray*}
It follows that there exist $c_4, c_5, c_6>0$ such that 
$$
\int_0^T V(t,u)dt -\int_0^T |W(t,u)|dt \geq c_4 |\bar{u}|^{\mu} - c_5 |\bar{u}|^{\nu} - c_6,
$$
and using that $\mu>\nu$ we deduce that $({\cal A}_{\infty})$ is satisfied. 

The proof that $(VW_\infty^1)$  implies $({\cal A}_{\infty})$ follows analogously.
\end{proof}

\begin{lemma}\label{VW1+A03}
Assume $(VW_\infty^1)$ and $({\cal A}_0^3)$ hold. Then condition $({\cal A}_0^1)$ is satisfied.
\end{lemma}

\begin{proof}
Since $(VW_\infty^1)$ holds, for any $(t,\xi)\in [0,T]\times \RR^3$ with $|\xi|\ge \rho$ one has 
$$
V(t,\xi)+|W(t,\xi)|\leq -\theta |\xi|^\mu+\theta |\xi|^\nu\to-\infty, \mbox{ as }|\xi|\to\infty,
$$
hence we can choose $R>\max\left\{r,1+\frac{T}{2},\rho\right\}$ such that 
$$
\max_{[0,T]\times \Gamma_R} (V+|W|)<\frac{1}{T}\min_{|\xi|=r}\int_0^T V(t,\xi)dt,
$$
which combined with $({\cal A}_0^3)$ shows that $({\cal A}_0^1)$ is satisfied.
\end{proof}

\begin{lemma}\label{ALPimplainf}
	If $(V, W)$ satisfies $(ALP),$ then $({\cal A}_{\infty})$ is fulfilled.
\end{lemma}
\begin{proof} Consider $c_1>0$ such that $|\nabla V|\leq c_1$ on $[0,T]\times\RR^3$ and assume 
	$$\int_0^TV(t,\xi)dt\to +\infty \  \ \mbox{as} \  \ |\xi|\to\infty.$$
Let $u\in D(\psi)$ be fixed, so $||u'||_Y\leq 1$ and $||\tilde{u}||_Y\le \frac{T}{4}<T.$ 
Note that 
\begin{eqnarray*}
	\frac{d}{ds} V(t,\bar{u}+s\tilde{u}(t)) = \nabla V(t,\bar{u}+s\tilde{u}(t)) \cdot\tilde{u}(t), \quad\mbox{for all} \ s\in [0,1], \ t\in [0,T],
\end{eqnarray*}
hence, for all $t\in [0,T],$
$$V(t,u(t)) - V(t,\bar{u}) = \int_0^1 \nabla V(t,\bar{u}+s\tilde{u}(t)) \cdot\tilde{u}(t)ds.$$
It follows that
\begin{eqnarray*}
	\int_0^T V(t,u)dt &=& \int_0^T V(t,\bar{u})dt + \int_0^T\int_0^1 \nabla V(t,\bar{u}+s\tilde{u}(t)) \cdot\tilde{u}(t)dsdt \\
	&\geq& \int_0^T V(t,\bar{u})dt - \int_0^T\int_0^1 |\nabla V(t,\bar{u}+s\tilde{u}(t))| |\tilde{u}(t)|dsdt \\
	&\geq& \int_0^T V(t,\bar{u})dt -T^2 c_1,
\end{eqnarray*}
and $({\cal A}_{\infty})$ follows using also that $|W|$ is bounded. The other situation follows analogously.
\end{proof}
Now the Corollaries follow from Theorem \ref{MainRes} and the above three Lemmas.

\section{Examples}

In this section we provide some explicit and nontrivial examples of potentials $V$ and $W$ that satisfy the two interactions required by the main result.

\begin{example}
Let $e,f\in\mathbb{R}^3$ with $|e|=|f|=1$, $\beta\neq  0$ and $\kappa>0$ such that 
$$
\kappa>\max\left\{ \frac{2(1+3T)}{T^4}+\frac{|\beta|}{2\pi T^3}, \frac{|\beta|}{T^2}\right\}
$$
and define the electromagnetic potential $(V,W)$ 
$$
V(t,\xi)=-\kappa\left(|\xi|^2 -T^2\right)^2- \cos\left(\frac{2\pi t}{T}\right) e\cdot \xi, \quad W(t,\xi)=\frac{\beta}{2}f\times \xi,
$$
for all $(t,\xi)\in [0,T]\times \RR^3$. Then $(VW_\infty^1)$ and $({\cal A}_0^3)$ hold, i.e., the assumptions of Corollary \ref{cor1} are satisfied. 
\end{example}

\begin{proof}
\begin{description}

\item {\sc Claim 1.} {\it $(VW_\infty^1)$ is satisfied with $\mu=4$, $\nu=3$, $\rho=2T$ and $\theta=\frac{|\beta|}{8T^2}$}.

Fix $t\in[0,T]$ and $\xi\in \RR^3$ with $|\xi|\ge \rho$. Then 
$$
|\xi|^2\ge 4T^2 \mbox{ and } |\xi|^3\ge 8T^3.
$$
Thus, 
$$
|W(t,\xi)|  \le \frac{|\beta|}{2}|f \times \xi|\leq \frac{|\beta|}{2}|\xi| \le \frac{|\beta|}{2}|\xi|\frac{|\xi|^2}{4T^2}=\theta |\xi|^{\nu}.
$$
On the other hand,
$$
|\xi|^2-T^2\ge |\xi|^2-\frac{|\xi|^2}{4}=\frac{3}{4}|\xi|^2, \mbox{ so } \left( |\xi|^2-T^2\right)^2\ge \frac{9}{16}|\xi|^4.
$$
Also,
$$
\kappa T^4\ge 2(1+3T)>6T, \mbox{ so }\frac{\kappa}{48}>\frac{1}{8T^3}.
$$
Consequently,
\begin{eqnarray*}
V(t,\xi) & \le & -\frac{9}{16}\kappa |\xi|^4 +\left| \cos\left( \frac{2\pi t}{T}\right)e\cdot \xi\right|
	\le -\frac{9}{16}\kappa |\xi|^4 +\frac{|\xi|^4}{|\xi|^3}\\
	&\le& -\frac{9}{16}\kappa |\xi|^4 +\frac{|\xi|^4}{8T^3} < -\frac{9}{16}\kappa |\xi|^4+\frac{\kappa}{48}|\xi|^4\\
	&=&-\frac{13}{24}\kappa |\xi|^4.
\end{eqnarray*}
Taking into account the choice of $\kappa$, one has 
$$
\theta=\frac{|\beta|}{8T^2}\le \frac{\kappa}{8}<\frac{13}{24}\kappa,
$$
hence $V(t,\xi)\leq -\theta \|\xi\|^\mu,$ and the Claim is now proved.

\item {\sc Claim 2.} {\it $({\cal A}_0^3)$ is satisfied with $r=T$.}

For any $(t,\xi) \in [0,T]\times\Gamma_T$ one has
$$
\left( |\xi|^2-T^2\right)^2\ge\frac{225}{256}T^4, \ \left|\cos\left( \frac{2\pi t}{T}\right)e\cdot \xi\right|\le |\xi| \ \mbox{ and } |W(t,\xi)|\leq \frac{|\beta|}{2}|\xi|.
$$
Now, the choice of $\kappa$ ensures that 
$$
\frac{225}{256} \kappa T^4> \frac{225}{256}\left(2+6T+\frac{|\beta|T}{2\pi}\right)>\frac{225}{128}+\frac{T}{4}+\frac{|\beta|T}{8}.
$$
Therefore, on $[0,T]\times \Gamma_T$
$$
V+|W|\le -\kappa \frac{225}{256}T^4+\frac{T}{4}+ \frac{|\beta|}{2}\frac{T}{4}<-\frac{225}{128}
$$
hence 
$$
T \max_{[0,T]\times \Gamma_T}\left( V+|W|\right)\le -\frac{225}{128}T<0.
$$
On the other hand, for any $|\xi|=r$ we have
\begin{eqnarray*}
\int_0^T V(t,\xi)dt & = &-\kappa \left(|\xi|^2-T^2 \right)^2 T- e\cdot \xi \int_0^T \cos\left(\frac{2\pi t}{T} \right)dt=0,
\end{eqnarray*}
and the Claim follows.

\end{description}
\end{proof}

\medskip

\begin{example}
Let $e,f\in\RR^3$ with $|e|=|f|=1$, $\beta\neq 0$ and  $\kappa>0$ such that 
$$
\kappa>\frac{1+T}{T^6}+\frac{|\beta|}{4\pi T^5}
$$
and consider the the electromagnetic potential 
$$
V(t,\xi)=\kappa \left( |\xi|^2-T^2\right)^2\left( |\xi|^2-\tau^2\right)-\cos\left( \frac{2\pi t}{T}\right)e\cdot \xi \mbox{ and }
W(t,\xi)=\frac{\beta}{2}f\times \xi,
$$
with $\tau=1+3T$.

Then $(VW_\infty^2)$ and $({\cal A}_0^1)$ are fulfilled, i.e., the assumptions of Corollary \ref{cor2} hold. 
\end{example}

\begin{proof}
\begin{description}
\item {\sc Claim 1}. {\it $(VW_\infty^2)$ is satisfied with $\mu=6$, $\nu=5$, $\theta=\frac{\kappa}{2}$ and $\rho=\max\left\{ 3\tau,\left(\frac{|\beta|}{\kappa}\right)^{1/4}\right\}$}.

Fix $(t,\xi)\in [0,T]\times \RR^3$ with $|\xi|\geq \rho$. Then 
$$
|\xi|\geq 3\tau \mbox{ and }|\xi|^4\geq \frac{|\beta|}{\kappa}
$$  
and 
$$
|W(t,\xi)|\le \frac{|\beta|}{2}|\xi| \le \frac{\kappa}{2}\frac{|\beta|}{\kappa}|\xi|=\frac{\kappa}{2}|\xi|^4|\xi|=\theta |\xi|^\nu.
$$
For $V$ we have the following estimates
$$
|\xi|^2-T^2\ge |\xi|^2-\tau^2\ge \frac{8}{9}|\xi|^2,
$$
and 
$$
(|\xi|^2-T^2)^2(|\xi|^2-\tau^2)\ge \frac{2}{3}|\xi|^6.
$$
The choice of $\kappa$ implies that 
$$
\kappa T^7>(1+T)T, \mbox{ so }\frac{1}{\kappa}<T^5 \mbox { and }|\xi|^5>(1+3T)^5>32T^5>\frac{32}{\kappa}.
$$  
Thus,
$$
\left| \cos\left(\frac{2\pi t}{T} \right)e\cdot \xi\right|\le |\xi|=\frac{|\xi|^6}{|\xi|^5}<\frac{\kappa}{32}|\xi|^6.
$$
Finally, 
$$
V(t,\xi) \ge \frac{2\kappa}{3}|\xi|^6-\frac{\kappa}{32}|\xi|^6=\frac{61}{96}\kappa |\xi|^6\geq \frac{\kappa}{2}|\xi|^6=\theta |\xi|^\mu.
$$

\item {\sc Claim 2.} {\it $({\cal A}_0^1)$  holds with $r=T$ and $R=1+\frac{5T}{2}$.}

By the choice of $\kappa$ 
$$
\kappa T^6>1+T+\frac{|\beta|}{4\pi} T \mbox{ and } \kappa T^5>1+\frac{|\beta|}{4\pi}.
$$
Direct computations yield 
$$
\min_{|\xi|=r}\int_0^T V(t,\xi)dr =0.
$$
Now, pick $(t,\xi)\in [0,T]\times \Gamma_T$. Then by the choice of $\tau$
$$
|\xi|^2-\tau^2\le -8T^2 \mbox{ and }\left(|\xi|^2-T^2\right)^2\ge \frac{225}{256}T^4.
$$
Thus
$$
V(t,\xi)\le -8\kappa \frac{225}{256}T^6+|\xi|\le -\frac{225}{32}T^6+\frac{T}{4}.
$$
Consequently, on $[0,T]\times \Gamma_T$
$$
V+|W| \leq  -\frac{225}{32}kT^6+\frac{T}{4}+\frac{|\beta|}{8}T\leq -\frac{225}{32}\left( 1+T+\frac{|\beta|}{4\pi}T\right)+\frac{T}{4}+\frac{|\beta|}{8}T<-\frac{225}{32}.
$$
Fix now $(t,\xi)\in [0,T]\times \Gamma_R$, i.e., $2T\le |\xi|\le 1+\frac{5T}{2}$. Then 
$$
\left(|\xi|^2-T^2 \right)^2 \ge 9T^4 \mbox{ and } \tau^2-|\xi|^2\ge\tau^2-\left(1+\frac{5T}{2}\right)^2 = T+\frac{11}{4}T^2.
$$
Therefore, on $[0,T]\times \Gamma_R$ we have
\begin{eqnarray*}
V+|W| & \le & -\kappa T^5 \left(9+\frac{99}{4}T\right)+|\xi|+\frac{|\beta|}{2}|\xi|\\
	&\le& -\left(1+\frac{|\beta|}{4\pi} \right)\left(9+\frac{99}{4}T \right)+\left( 1+\frac{|\beta|}{2}\right)\left( 1+\frac{5T}{2}\right)\\
	& \le & -8.
\end{eqnarray*}
Finally, 
$$
T\max_{[0,T]\times (\Gamma_T\cup \Gamma_R)}(V+|W|)<0=\min_{|\xi|=r}\int_0^T V(t,\xi)dt.
$$
\end{description}
\end{proof}

\medskip

The only thing missing now is a proof for the Minimax theorem for the Poincar\'{e} action functional, namely Theorem \ref{MinMaxThm}, which is provided in the following two sections: first we prove a deformation result which requires only condition $(C)$, then we use the deformation together with property $(S)$ to close the proof.

\section{Deformations for the Poincar\'{e} action}\label{SecDeformation}

Throughout this section we combine ideas from Alves \& de Morais Filho \cite{AlvFil}, Arcoya, Bereanu \& Torres \cite{ABT,ABT2}, Byeon, Ikoma, Malchiodi \& Mari \cite{BIMM}, Schechter \& Zou \cite{Schechter-Zou-JFA03} and Szulkin \cite{Szulkin} to prove a deformation theorem, with localization information, designed for linking sets that are not compact, but merely closed and bounded. We also point out the fact that condition $(C)$ is the only compactness we need.

We start with some notation used in the remainder of the paper and collect some elementary consequences of the qualitative properties of our functional that are used repeatedly. 

For $Q\subset X$ and $\delta>0$ we write
$$
[Q]^Y_\delta=\left\{u\in X:\dist_Y(u,Q)\le\delta\right\}
$$
for the closed $\delta$-neighbourhood of $Q$ in the $Y$-metric.  The set of all critical points of $\cal F$ is denoted by  ${\cal K}_{\F}$ and we write 
$$
{\cal K}_{\F}(c)={\cal K}\cap \left\{u\in X: \F(u)=c\right\},
$$
for a given value $c\in\mathbb{R}$.

For $u\in D(\psi)$ and $v\in X$ we define the functional
\begin{equation}\label{DescFun}
\D(u;v)=\langle\phi'(u),v-u\rangle+\psi(v)-\psi(u).
\end{equation}
We use the convention $\D(u;v)=+\infty$ when $v\notin D(\psi)$ (then $\psi(v)=+\infty$), so that $\D(u;\cdot)$ is proper and convex. Symmetrically, we set $\D(u;v)=-\infty$ whenever $u\notin D(\psi)\mbox{ and }v\in D(\psi).$ We also point out the fact that the case $u,v\not\in D(\psi)$ never occurs throughout the paper.

Let $D(\partial\psi)=\{u\in X:\partial\psi(u)\neq\varnothing\}$. It is well known that for a convex and lower semicontinuous functional $D(\partial\psi)\subseteq D(\psi)$ and $\operatorname{int}D(\partial\psi)=\operatorname{int}D(\psi)$. Define $\lF:X\to[0,+\infty]$ by
$$
\lF(u)=
\begin{cases}
\displaystyle\inf_{\xi\in\partial\psi(u)}\|\phi'(u)+\xi\|_{X^*}, & u\in D(\partial\psi),\\[2mm]
+\infty, & u\notin D(\partial\psi).
\end{cases}
$$
Thus $\lF(u)$ is the distance in $X^*$ from $-\phi'(u)$ to the convex set $\partial\psi(u)$, and plays the role of the norm of the derivative of $\F$ at $u$.

\begin{definition}\label{PSc-seq}
A sequence $\{u_n\}\subset D(\psi)$ is a {\it Palais--Smale sequence at level $c\in\mathbb R$} for $\F$, abbreviated $(PS)_c$-sequence, if
\begin{description}
\item $\bullet$ $\F(u_n)\to c$ as $n\to \infty$;

\item  $\bullet$ there exists $\{\varepsilon_n\}\subset (0,\infty)$ such that  $\varepsilon_n\searrow0$ and
\begin{equation}\label{PScond}
\langle\phi'(u_n),v-u_n\rangle+\psi(v)-\psi(u_n)\ge-\varepsilon_n\|v-u_n\|_X,\  \forall v\in X,\ \forall n\geq 1.
\end{equation}
\end{description}
\end{definition}

\begin{lemma}\label{LambdaProp}
Let $u\in D(\psi)$. Then the following hold.
\begin{enumerate}[$(a)$]
\item For every $v\in X$,
$$
\D(u;v)\ge-\lF(u)\,\|v-u\|_X,
$$
and, moreover,
\begin{equation}\label{SlopeId}
\lF(u)=\sup_{v\in X\setminus\{u\}}\frac{\big(-\D(u;v)\big)_+}{\|v-u\|_X},\  (x)_+=\max\{x,0\}.
\end{equation}
Consequently, for every $\mu\ge0$,
\begin{equation}\label{SlopeEquiv}
\lF(u)\le\mu\iff \D(u;v)\ge-\mu\|v-u\|_X\ \ \forall v\in X.
\end{equation}

\item The following assertions are equivalent:
\begin{description}
\item (i) $u$ is a critical point of $\F$;
\item (ii) $\D(u;v)\ge0$ for all $v\in X$;
\item (iii) $-\phi'(u)\in\partial\psi(u)$;
\item (iv) $\lF(u)=0$.
\end{description}

\item For a sequence $\{u_n\}\subset D(\psi)$ the following statements are equivalent:
\begin{description}
\item (i') $\{u_n\}$ is a $(PS)_c$-sequence for $\F$;
\item (ii') $\F(u_n)\to c$ and there exists a sequence of positive numbers $\{\varepsilon_n\}$ such that  $\varepsilon_n\to0$ and $\D(u_n;v)\ge-\varepsilon_n\|v-u_n\|_X$, for all $v\in X$ and all $n\geq 1$;
\item (iii') $\F(u_n)\to c$ and $\lF(u_n)\to0$ as $n\to\infty$.
\end{description}
\end{enumerate}
\end{lemma}

\begin{proof}
Fix $u\in D(\psi)$.
\begin{description}
\item [$(a)$]  If $u\notin D(\partial\psi)$ then $\lF(u)=+\infty$ and the inequality is trivial. If $u\in D(\partial\psi)$, then fix $\xi\in\partial\psi(u)$. Since $\psi(v)-\psi(u)\ge\langle\xi,v-u\rangle$ for all $v\in X$ one has
\begin{eqnarray*}
\D(u;v) &=& \langle\phi'(u),v-u\rangle+\psi(v)-\psi(u)\ge \langle\phi'(u)+\xi,\,v-u\rangle\\ 
	&\ge& -\|\phi'(u)+\xi\|_{X^\ast}\|v-u\|_X.
\end{eqnarray*}
Taking the infimum over $\xi\in\partial\psi(u)$ yields $\D(u;v)\ge-\lF(u)\|v-u\|_X$.

For the identity \eqref{SlopeId} we write $s_\F(u)$ for the right-hand side of the equality. From the inequality we just proved, $-\D(u;v)\le\lF(u)\|v-u\|_X$, hence
$$
s_\F(u)\ \le\ \lF(u).
$$
For the reverse inequality assume $s_\F(u)<+\infty$ (otherwise there is nothing to prove). By the definition of $s_\F$,
$$
\D(u;v)\ \ge\ -s_\F(u)\|v-u\|_X, \quad\mbox{for all }v\in X,
$$
that is, the convex function
$$
g(v)=\psi(v)+\langle\phi'(u),v-u\rangle+s_\F(u)\|v-u\|_V
$$
satisfies $g(v)\ge\psi(u)=g(u)$ for all $v$, therefore $u$ minimizes $g$ and $0\in\partial g(u)$. The last two summands of $g$ are convex and continuous on $X$, so the Moreau--Rockafellar sum rule applies and gives
$$
0\in \partial\psi(u)+\phi'(u)+s_\F(u) \partial\left(\|\cdot-u\|_X\right)(u)\ =\ \partial\psi(u)+\phi'(u)+s_\F(u) B_{X^*}[0,1],
$$
because the subdifferential of the norm at the origin is the closed dual unit ball $B_{X^*}[0,1]$. Thus there exist $\xi\in\partial\psi(u)$ and $b\in B_{X^*}[0,1]$ with $\phi'(u)+\xi+s_\F(u)b=0$, whence
$$
\lF(u)\le \|\phi'(u)+\xi\|_{X^\ast} =s_\F(u)\|b\|_{X^\ast} \le s_\F(u).
$$
This proves \eqref{SlopeId}. The equivalence \eqref{SlopeEquiv} is a direct reformulation, since $\lF(u)\le\mu$ says precisely that every quotient in \eqref{SlopeId} is less than or equal to $\mu$.

\item [$(b)$] The equivalence $(i)\Leftrightarrow(ii)$ is exactly the definition of a critical point  (see Theorem \ref{LFEPoincare}). On the other hand,  $(i)\Leftrightarrow(iii)$ is the characterisation $0\in\phi'(u)+\partial\psi(u)$ of a Szulkin critical point \cite{Szulkin}. It remains to show $(iii)\Leftrightarrow(iv)$.

If $(iii)$ holds, then $u\in D(\partial\psi)$ and, taking $\xi=-\phi'(u)$ in the definition of $\lF$,
$$
\lF(u)\ \le\ \|\phi'(u)+(-\phi'(u))\|_{X^\ast} = 0.
$$
Conversely, assume $(iv)$. Then $u\in D(\partial\psi)$ and there exists a sequence $\{\xi_n\}\subset\partial\psi(u)$ such that $\|\phi'(u)+\xi_n\|_{X^\ast}\to0$, i.e., $\xi_n\to-\phi'(u)$ in $X^\ast$. Since $\partial\psi(u)$ is convex and weak$^\ast$-closed, hence norm-closed, the limit $-\phi'(u)$ lies in $\partial\psi(u)$, which is $(iii)$.

\item [$(c)$] $(i')\Leftrightarrow(ii')$ is Definition~\ref{PSc-seq}.

\smallskip
$(ii')\Rightarrow(iii')$. Fix $n$. By hypothesis $\D(u_n;v)\ge-\varepsilon_n\|v-u_n\|_X$ for all $v$, so \eqref{SlopeEquiv} (with $\mu=\varepsilon_n$) gives $\lF(u_n)\le\varepsilon_n$. Hence $0\le\lF(u_n)\le\varepsilon_n\to0$.

\smallskip
$(iii')\Rightarrow(ii')$. Set
$$
\varepsilon_n=\sup_{k\ge n}\lF(u_k)+\tfrac1n .
$$
Since $\lF(u_k)\to0$, the supremum is finite, and $\{\varepsilon_n\}$ is strictly positive, non-increasing and tends to $0$, with $\varepsilon_n\ge\lF(u_n)$. By part $(a)$,
$$
\D(u_n;v)\ge-\lF(u_n)\|v-u_n\|_X\ge-\varepsilon_n\|v-u_n\|_X, \quad\mbox{for all }v\in X,\ n\ge1,
$$
which is $(ii')$.
\end{description}
\end{proof}

\begin{lemma}\label{Ysemicont}
Let $\{u_n\}\subset D(\psi)$ be bounded with $u_n\to u$ in $Y$. The following assertions hold:
\begin{enumerate}[$(i)$]

\item $u\in D(\psi)$;

\item $\F(u)\le\liminf\limits_{n\to \infty}\F(u_n)$;

\item $\limsup\limits_{n\to\infty}\D(u_n;v)\le\D(u;v)$, for every $v\in X$.
\end{enumerate}
\end{lemma}

\begin{proof}
\begin{description}

\item $(i)$ Using Lemma \ref{CompEmbed} we infer that $u\in X$. Since $\{u_n\}\subset D(\psi)$, property $(P_5)$ gives
$$
\psi(u)\leq \liminf_{n\to \infty}\psi(u_n)\leq T<+\infty,
$$
hence $u\in D(\psi)$.

\item $(ii)$ By property $(P_2)$ we have $\phi(u_n)\to \phi(u)$. This together with $(P_5)$ gives
\begin{align*}
\liminf_{n\to\infty} \F(u_n) & = \liminf_{n\to\infty} \left(\phi(u_n)+\psi(u_n) \right)\\
	& \geq \liminf_{n\to\infty} \phi(u_n)+\liminf_{n\to\infty} \psi(u_n)\\
	& \geq \phi(u)+\psi(u)=\F(u).
\end{align*}

\item $(iii)$ Let $v\in X$ be fixed. If $v\in X\setminus D(\psi)$, then
$$
\D(u;v)=\langle \phi'(u),v-u\rangle+\psi(v)-\psi(u)=+\infty
$$
and the conclusion follows immediately.

Assume now $v\in D(\psi)$. Then
\begin{align*}
\D(u_n;v) & = \langle \phi'(u_n),v-u_n\rangle+\psi(v)-\psi(u_n) \\
	& = \langle \phi'(u_n),v\rangle -\langle \phi'(u_n),u_n\rangle+\psi(v)-\psi(u_n).
\end{align*}
By $(P_2)$ and $(P_5)$ we have
$$
\lim_{n\to \infty}\langle \phi'(u_n),v \rangle = \langle \phi'(u),v\rangle, \ \lim_{n\to \infty}\langle \phi'(u_n),u_n \rangle = \langle \phi'(u),u\rangle
$$
and
$$
\liminf_{n\to\infty} \psi(u_n)\geq \psi(u).
$$
Therefore,
\begin{align*}
\limsup_{n\to\infty}\D(u_n;v) & = \limsup_{n\to\infty}\left( \langle \phi'(u_n),v\rangle -\langle \phi'(u_n),u_n\rangle+\psi(v)-\psi(u_n)\right)\\
	&\leq \limsup_{n\to\infty} \langle \phi'(u_n),v\rangle -\liminf_{n\to\infty} \langle \phi'(u_n),u_n\rangle+\psi(v)-\liminf_{n\to\infty}\psi(u_n)\\
	& \leq \langle \phi'(u),v\rangle-\langle \phi'(u),u\rangle+\psi(v)-\psi(u)\\
	&=\D(u;v).
\end{align*}
\end{description}
\end{proof}

\begin{lemma}\label{BandBounded}
Assume condition $(C)$ holds. Then, for every $c\in\mathbb R$ and $\beta>0$, the strip
$$
\Sigma_{c,\beta}=\{u\in D(\psi):|\F(u)-c|\le3\beta\}
$$
is bounded in $X$. Moreover, its closure $\overline{\Sigma_{c,\beta}}^{\,Y}$ in $(Y,\|\cdot\|_Y)$ is compact and contained in $\{u\in X:\|u\|_X\le R_{c,\beta}\}$, where $R_{c,\beta}$ is any bound for $\Sigma_{c,\beta}$.
\end{lemma}

\begin{proof}
Using property $(P_5)$, on $\Sigma_{c,\beta}$ one has 
$$
c-3\beta-T\leq \phi=\F-\psi\leq c+3\beta,
$$
so $|\phi|$ is bounded there. 

Arguing be contradiction, suppose $\Sigma_{c,\beta}$ is unbounded. Then  there exists $\{u_n\}\subset \Sigma_{c,\beta}$ with $\|u_n\|_X\to+\infty$. 
By Lemma \ref{basicineq} one has 
$$
|u_n(t)|\leq |\bar u_n| +|\tilde u_n(t)|\le |\bar u_n|+\|\tilde u_n\|_Y\le |\bar u_n|+\frac{T}{4}.
$$
Taking the supremum as $t\in [0,T]$ and using again Lemma \ref{basicineq} we conclude that 
$$
|\bar u_n|\geq \|\tilde u_n\|_Y-\frac{T}{4}\geq \|u_n\|_X-1-\frac{T}{4}\to +\infty, \mbox{ as }n\to+\infty.
$$
Now, by condition $(C)$  we have $|\phi(u_n)|\to+\infty$, which is impossible on $\Sigma_{c,\beta}$. 
Thus $\Sigma_{c,\beta}$ is bounded, say $\Sigma_{c,\beta}\subset B_X[0,R_{c,\beta}]$. By the compact embedding $X\hookrightarrow Y$ this set is relatively compact in $Y$, so $\overline{\Sigma_{c,\beta}}^{\,Y}$ is $Y$-compact. Since closed balls of $X$ are $Y$-closed due to Lemma \ref{CompEmbed}, it follows that $\overline{\Sigma_{c,\beta}}^{\,Y}\subset B_X\left[0,R_{c,\beta}\right]$.
\end{proof}

The next lemma turns a $(PS)_c$-sequence into a genuine critical point \textit{at the level $c$}: the boundedness of the sequence is supplied by Lemma~\ref{BandBounded}, the $Y$-limit by the compact embedding $X\hookrightarrow Y$, the criticality by Lemma~\ref{Ysemicont}$(iii)$, and the identification of the level by the third convergence in $(P_2)$. 

\begin{lemma}\label{PScritical}
Assume condition $(C)$ holds, fix $c\in\mathbb R$ and let $\{u_n\}\subset D(\psi)$ be a $(PS)_c$-sequence for $\F$. Then:
\begin{description}
\item (i) the sequence $\{u_n\}$ is bounded in $X$, and there exist a subsequence $\{u_{k_n}\}$ and $u^\ast\in X$ such that $\|u_{k_n}-u^\ast\|_Y\to0$;
\item (ii) for every such $Y$-limit $u^\ast$ one has $u^\ast\in\mathcal K_{\cal F}$;
\item (iii) moreover $\psi(u_{k_n})\to\psi(u^\ast)$, $\phi(u_{k_n})\to\phi(u^\ast)$ and consequently
$$
\F(u^\ast)=\lim_{n\to\infty}\F(u_{k_n})=c .
$$
\end{description}
\end{lemma}

\begin{proof}
\begin{description}
\item $(i)$ Since $\F(u_n)\to c$, there exists $n_0\geq 1$ such that $|\F(u_n)-c|\le3\beta$ for all $n\ge n_0$ (say with $\beta=1$), i.e., $u_n\in\Sigma_{c,1}$ for $n\ge n_0$. By Lemma~\ref{BandBounded} the sequence $\{u_n\}_{n\ge n_0}$, hence $\{u_n\}$ itself, is bounded in $X$. By the compact embedding $X\hookrightarrow Y$, a subsequence $\{u_{k_n}\}$ converges in $Y$ to some $u^\ast$, and $u^\ast\in X$  due to Lemma \ref{CompEmbed}.

\item $(ii)$  Relabel $\{u_{k_n}\}$ as $\{u_n\}$. By Lemma~\ref{Ysemicont}$(i)$, $u^\ast\in D(\psi)$. Fix $v\in X$. Since $\{u_n\}$ is bounded, $\|v-u_n\|_X\le\|v\|_X+\sup_n\|u_n\|_X=:C_v<+\infty$, so \eqref{PScond} gives $\D(u_n;v)\ge-\varepsilon_n C_v\to0$, whence $\liminf_n\D(u_n;v)\ge0$. By Lemma~\ref{Ysemicont}$(iii)$,
$$
\D(u^\ast;v)\ \ge\ \limsup_n\D(u_n;v)\ \ge\ \liminf_n\D(u_n;v)\ \ge\ 0 .
$$
As $v\in X$ was arbitrary, $u^\ast\in\mathcal K$.

\item $(iii)$ Taking $v=u^\ast\in D(\psi)$ in \eqref{PScond} we have 
\begin{equation}\label{LevelTest}
\psi(u_n)\ \le\ \psi(u^\ast)+\langle\phi'(u_n),u^\ast-u_n\rangle+\varepsilon_n\|u^\ast-u_n\|_X,\forall n\ge 1.
\end{equation}
The sequence $\{u_n\}$ being bounded and $\varepsilon_n\to0$, the last term tends to $0$. For the middle one, $(P_2)$ gives $\phi'(u_n)\to\phi'(u^\ast)$ in $X^*$, hence $\langle\phi'(u_n),u^\ast\rangle\to\langle\phi'(u^\ast),u^\ast\rangle$, and also $\langle\phi'(u_n),u_n\rangle\to\langle\phi'(u^\ast),u^\ast\rangle$. Therefore
$$
\langle\phi'(u_n),u^\ast-u_n\rangle\to0 .
$$
Passing to the $\limsup$ in \eqref{LevelTest} we obtain $\limsup_n\psi(u_n)\le\psi(u^\ast)$, which combined with the lower semicontinuity in $(P_5)$, $\liminf_n\psi(u_n)\ge\psi(u^\ast)$, yields $\psi(u_n)\to\psi(u^\ast)$. Since $\phi(u_n)\to\phi(u^\ast)$ by $(P_2)$, we get 
$$
\F(u_n)=\phi(u_n)+\psi(u_n)\to\phi(u^\ast)+\psi(u^\ast)=\F(u^\ast),
$$
and $\F(u_n)\to c$ forces $\F(u^\ast)=c$.
\end{description}
\end{proof}

\begin{remark}\label{NoPS}
Lemma~\ref{PScritical} shows that under condition $(C)$ the functional $\F$ automatically satisfies the following weak Palais--Smale condition at every level $c\in\mathbb R$, which was introduced in \cite{ABT2}:

\noindent {\bf ${\bf(wPS)_c}$-condition:} {\it  every $(PS)_c$-sequence possesses a subsequence which converges in $Y$ to a critical point of $\F$ of level $c$. }
\end{remark}

\begin{lemma}\label{Equiv}
Let $c\in\mathbb R$ be fixed and $ Q \subset X$. The following assertions are equivalent:
\begin{description}
\item$(i)$ $\F$ possesses no $(PS)_c$-sequence satisfying  $\dist_Y(u_n, Q )\to0$ as $n\to\infty$;

\item$(ii)$ there exist $\beta\in(0,1)$, $\delta_0\in(0,\beta)$ and $\rho>0$ such that for every $u\in D(\psi)$ satisfying
$$
|\F(u)-c|\leq 3\beta \mbox{ and } \dist_Y(u, Q )\leq \rho
$$
there exists $v=v(u)$ with the property that
\begin{equation}\label{SlopeBand}
\D(u;v)<-5\delta_0\|v-u\|_X,
\end{equation}
equivalently $(\mbox{due to }$\eqref{SlopeEquiv}$)$ $\lF(u)>5\delta_0$.
\end{description}
\end{lemma}

\begin{proof}
\begin{description}

\item$(i)\Rightarrow (ii).$ Assume by contradiction that \eqref{SlopeBand} fails for every choice of parameters.
Then, for each positive integer $n\in\mathbb{N}$,  choosing $\beta=\frac{1}{3n}$, $\delta_0=\frac{1}{5n}$ and $\rho=\frac{1}{n}$ produces  $u_n\in D(\psi)$ with
$$
 |\F(u_n)-c|\le\frac{1}{n}, \ \dist_Y(u_n, Q )\leq \frac{1}{n} \mbox{ and }\D(u_n;v)\geq-\frac{1}{n}\|v-u_n\|_X,\quad \forall v\in X,
$$
i.e., a $(PS)_c$-sequence with $\dist_Y(u_n, Q )\to 0$, contradicting $(i)$.

\item$(ii)\Rightarrow (i).$ Assume, again by contradiction, that $\{u_n\}\subset D(\psi)$ is a $(PS)_c$-sequence such that $\dist_Y(u_n, Q )\to 0$. Then $\F(u_n)\to c$ and there exists $\{\varepsilon_n\}\subset (0,\infty)$ with $\varepsilon_n\searrow 0 $ and $\D(u_n;w)\geq -\varepsilon_n\|w-u_n\|_X$ for all $w\in X$. Take  $n_0\in\mathbb{N}$ sufficiently large so that
$$
|\F(u_n)-c|\le3\beta,\ \dist_Y(u_n, Q )\leq \rho \mbox{ and }\varepsilon_n\leq 5\delta_0,\quad \forall n\geq n_0.
$$
For $n\ge n_0$, apply $(ii)$ to $u=u_n$ and let $v(u_n)$ be the element provided by \eqref{SlopeBand}. Then, on one hand 
$$
\D(u_n;v(u_n))\ge -\varepsilon_n\|v(u_n)-u_n\|_X\ge -5\delta_0\|v(u_n)-u_n\|_X,
$$
while on the other hand \eqref{SlopeBand} gives $\D(u_n;v(u_n))<-5\delta_0\|v(u_n)-u_n\|_X$, a contradiction.
\end{description}
\end{proof}

\begin{lemma}\label{ConvComb}
Fix $\Theta>0$. Let $w\in D(\psi)\cap B_X[0,\Theta]$, $N\in \mathbb{N}$, $v_j\in D(\psi)$ and  $\theta_j\ge0$, $1\leq j\leq N$, be such that  $\lambda=\sum_{j=1}^{N}\theta_j\le\tfrac12$. Set
$$
\eta^+=(1-\lambda)w+\sum_{j=1}^{N}\theta_j v_j,\quad \eta^-=(1+\lambda)w-\sum_{j=1}^{N}\theta_j v_j .
$$
Then $\eta^+\in D(\psi)$,
\begin{equation}\label{Forward}
\F(\eta^+)-\F(w)\le\sum_{j=1}^{N}\theta_j\,\D(w;v_j)+\mathcal{R}^+,
\end{equation}
and
\begin{equation}\label{Backward}
\F(\eta^-)-\F(w)\ge-\sum_{j=1}^{N}\theta_j\,\D(w;v_j)+\mathcal{R}^-,
\end{equation}
where 
$$
\mathcal{R}^+=\phi(\eta^+)-\phi(w)-\sum_{j=1}^{N}\theta_j\langle\phi'(w),v_j-w\rangle
$$ 
and 
$$
\mathcal{R}^-=\phi(\eta^-)-\phi(w)+\sum_{j=1}^{N} \theta_j \langle \phi'(w),v_j-w\rangle.
$$ 
Moreover, for every $\varepsilon>0$ there exists $r_\varepsilon=r_\varepsilon(\Theta)\in(0,1]$, depending only on $\Theta$ and $\varepsilon$, such that whenever in addition $\sum_{j=1}^{N}\theta_j\|v_j-w\|_X\le r_\varepsilon$ one has
\begin{equation}\label{R-Estimate}
|\mathcal{R}^{\pm}|\le\varepsilon\sum_{j=1}^{N}\theta_j\|v_j-w\|_X.
\end{equation}
\end{lemma}

\begin{proof}
Since $\lambda\le1$ and $w,v_j\in D(\psi)$, $\eta^+$ is a convex combination of points of the convex set $D(\psi)$, so $\eta^+\in D(\psi)$, and Jensen's inequality gives
$$
\psi(\eta^+)\leq(1-\lambda)\psi(w)+\sum_{j=1}^{N}\theta_j\psi(v_j),
$$
i.e.,
\begin{align*}
\psi(\eta^+)-\psi(w) & \leq -\lambda\psi(w)+\sum_{j=1}^{N}\theta_j\,\psi(v_j)=-\psi(w)\sum_{j=1}^{N} \theta_j +\sum_{j=1}^{N} \theta_j \psi(v_j)\\
		& = \sum_{j=1}^{N} \theta_j \left(\psi(v_j)-\psi(w) \right).
\end{align*}
Adding $\phi(\eta^+)-\phi(w)$ we get
\begin{align*}
\F(\eta^+)-\F(w) & \leq \sum_{j=1}^{N}\theta_j\,(\psi(v_j)-\psi(w)) +\phi(\eta^+)-\phi(w)\\
	& = \sum_{j=1}^{N}\theta_j\,(\psi(v_j)-\psi(w)) +\sum_{j=1}^{N} \theta_j \langle \phi'(w),v_j-w\rangle+\mathcal{R}^+\\
	& = \sum_{j=1}^{N} \theta_j \D(w;v_j)+\mathcal{R}^+,
\end{align*}
showing  \eqref{Forward} holds.

For \eqref{Backward}, note that $w=\tfrac{1}{1+\lambda}\eta^-+\sum_{j=1}^{N}\tfrac{\theta_j}{1+\lambda}v_j$ is a convex combination, so
$$
\psi(w)\le\tfrac1{1+\lambda}\psi(\eta^-)+\sum_{j=1}^{N}\tfrac{\theta_j}{1+\lambda}\psi(v_j),
$$
whence, multiplying by $1+\lambda$ and rearranging (this also covers the case $\eta^-\notin D(\psi)$, where the right-hand side is $+\infty$),
$$
\psi(\eta^-)-\psi(w)\ge-\sum_{j=1}^{N}\theta_j(\psi(v_j)-\psi(w)).
$$
Adding $\phi(\eta^-)-\phi(w)$ and using the definition of $\mathcal{R}^-$ gives
$$
\F(\eta^-)-\F(w)\ge-\sum_{j=1}^{N}\theta_j\D(w;v_j)+\mathcal{R}^-,
$$
which is \eqref{Backward}.

Now we prove \eqref{R-Estimate}. 
Assume $\sum_{j=1}^{N} \theta_j\|v_j-w\|_X\le r$ with $r\le1$ to be fixed later. Then $\|\eta^+-w\|_X\le\sum_{j=1}^{N}\theta_j\|v_j-w\|_X\le r\le1$, so the segment $[w,\eta^+]\subset B_X\left[0,\Theta+1\right]$. Applying the mean value theorem to $g:[0,1]\to\mathbb{R}$,
$$
g(t)=\phi(w+t(\eta^+-w)),
$$
we get $t^*\in(0,1)$ with $\phi(\eta^+)-\phi(w)=g(1)-g(0)=g'(t^*)=\langle\phi'(\omega),\eta^+-w\rangle$, where $\omega=w+t^*(\eta^+-w)\in B_X[0,\Theta+1]$. Since $\eta^+-w=\sum_{j=1}^{N}\theta_j(v_j-w)$,
$$
\mathcal{R}^+=\sum_{j=1}^{N}\theta_j\langle\phi'(\omega)-\phi'(w),v_j-w\rangle,
$$
and
$$
|\mathcal{R}^+|\le\|\phi'(\omega)-\phi'(w)\|_{X^*}\sum_{j=1}^{N}\theta_j\|v_j-w\|_X.
$$
By property $(P_4)$ applied on $B_X\left[0,\Theta+1\right]$, given $\varepsilon>0$ there exists $\delta>0$ such that $\|\phi'(\omega)-\phi'(w)\|_{X^*}\le\varepsilon$ whenever $\|\omega-w\|_Y\le\delta$. As 
$$
\|\omega-w\|_Y\le \|\omega-w\|_X\le \|\eta^+-w\|_X\le  r,
$$
choosing $r_\varepsilon=\min\{1,\delta\}$ yields 
$$
|\mathcal{R}^+|\le\varepsilon\sum_{j=1}^{N}\theta_j\|v_j-w\|_X.
$$
Note $r_\varepsilon$ depends only on $\Theta$ and $\varepsilon$. 

The bound for $\mathcal{R}^-$ is obtained in the same way, with $\eta^-$ in place of $\eta^+$ (here $\|\eta^--w\|_X=\|\sum_{j=1}^{N}\theta_j(v_j-w)\|_X\le r$ as well).
\end{proof}

We now fix the data used by the deformation. Let $Q\subseteq X$, let $c\in\mathbb R$ and $\beta\in(0,1)$, $\delta_0\in(0,\beta)$, $\rho>0$ be such that for every $u\in D(\psi)$ with $|\F(u)-c|\le3\beta$ and $\dist_Y(u,Q)\le\rho$ there exists $v(u)$ such that
\begin{equation}\label{SlopeStrip}
\D(u;v(u))<-5\delta_0\|v(u)-u\|_X.
\end{equation}
By Lemma~\ref{Equiv} this holds precisely when $\F$ has no $(PS)_c$-sequence approaching $Q$. The {\it global case} corresponds to $Q=X$ and an arbitrary $\rho>0$ (say $\rho=1$) and Lemma~\ref{Equiv} then says precisely that $\F$ possesses no $(PS)_c$-sequence at all. Note $v(u)\neq u$, since $\D(w;w)=0$.

Let $K\subset\{u\in D(\psi):|\F(u)-c|\le3\beta,\ \dist_Y(u,Q)\le\rho/2\}$ and $\widehat K$ be its closure in $Y$, and let $\widehat R>0$ bound the strip $\Sigma_{c,\beta}$ (such an $\widehat R$ is provided by Lemma~\ref{BandBounded}). Then $K\subset B_X\left[0,\widehat R\right]$. Since closed balls of $X$ are $Y$-closed (by Lemma \ref{CompEmbed}) and $B_X\left[0,\widehat R\right]$ is $Y$-precompact,
$$
\widehat K=\overline K^Y\subset B_X\left[0,\widehat R\right] \mbox{ is compact in }Y .
$$
By Lemma~\ref{Ysemicont}, every $w\in\widehat K$ lies in $D(\psi)$ with $\F(w)\le c+3\beta$ and every $w\in\widehat K$ satisfies $\dist_Y(w,Q)\le\rho/2$, so \eqref{SlopeStrip} is available at each such $w\in \widehat K$ with $\F(w)\ge c-3\beta$.
\begin{lemma}\label{DescSel}
Assume  \eqref{SlopeStrip} holds and fix $\Theta_0\ge \widehat R$. Then for each $u\in\widehat K$ there exist $\widehat v=\widehat v(u)\in D(\psi)$, $r=r(u)\in(0,1]$ and  $\alpha=\alpha(u)>0$ such that, for all $w\in X$ with $\|w\|_X\le\Theta_0$ and $\|w-u\|_Y<r$ one has:
\begin{description}

\item$(i)$ if $\F(u)\ge c-3\beta$, then $\D(w;\widehat v)\le-\alpha$, with $\alpha=\tfrac52\delta_0\|\widehat v-u\|_X>0$;

\item$(ii)$ if $\F(u)<c-3\beta$ and $u\notin\mathcal K_{\F}$, then $ \D(w;\widehat v)\le0$, and $\F(w)\ge c-2\beta\Rightarrow\D(w;\widehat v)\le-\alpha$, with  $\alpha=\tfrac{\beta}{4}$;

\item$(iii)$ if $\F(u)<c-3\beta$ and $u\in\mathcal K_{\F}$, then $ \D(w;\widehat v)\le\tfrac{\delta_0}{2}\|\widehat v-w\|_X$, and $\F(w)\ge c-2\beta\Rightarrow\D(w;\widehat v)\le-\alpha$, with $\alpha=\tfrac{\beta}{2}$.
\end{description}
In particular,  $\D(w;\widehat v)\le\tfrac{\delta_0}{2}\|\widehat v-w\|_X$ for every $w$ in the indicated set and in all cases 
$$
\F(w)\ge c-2\beta\Rightarrow\D(w;\widehat v)\le-\alpha.
$$ 
\end{lemma}

\begin{proof}
All points below lie in $B_X\left[0,\Theta_0\right]$, whose closure in $Y$ is compact. By Lemma~\ref{Ysemicont} and $\F(u)\le c+3\beta$ on $K$, every $u\in\widehat K$ has $\F(u)\le c+3\beta$.  If, in addition, $\F(u)\ge c-3\beta$, then $|\F(u)-c|\le3\beta$ and $\dist_Y(u,Q)\le\rho/2\le\rho$, so \eqref{SlopeStrip} applies at $u$.

\begin{description}
\item$(i)$ Here $|\F(u)-c|\le3\beta$ and $\dist_Y(u,Q)\le\rho/2$, so we may take $\widehat v=v(u)\ne u$ provided by \eqref{SlopeStrip}. Then $\widehat v\in D(\psi)$ and $\D(u; \widehat v)<-5\delta_0\|\widehat v-u\|_X=-2\alpha$ with $\alpha=\frac{5}{2}\delta_0\|\widehat v-u\|_X$.

If no radius worked, we could find $w_n\to u$ in $Y$ and $\|w_n\|_X\le\Theta_0$ such that $\D(w_n;\widehat v)>-\alpha$. But, by Lemma~\ref{Ysemicont}$(iii)$, $u\mapsto\D(u;\widehat v)$ is $Y$-upper semicontinuous along bounded sequences, hence
$$
-\alpha\le\limsup_{n\to\infty}\D(w_n;\widehat v)\le\D(u;\widehat v)<-2\alpha,
$$
a contradiction. Sequences with $w_n\notin D(\psi)$ may be discarded, since $\D(w_n;\widehat v)=-\infty$. Thus $\D(w;\widehat v)\le-\alpha$ for admissible $w$, and $-\alpha\le0\le\tfrac{\delta_0}{2}\|\widehat v-w\|_X$.

\item$(ii)$ Here $\F(u)<c-3\beta$ and $u\notin\mathcal K_{\F}$. Since $u\notin\mathcal K_{\F}$, there exists $\widehat v=v\in D(\psi)$, $v\neq u$, with $\D(u; \widehat v)=-2\kappa_u<0$.
As in case $(i)$, Lemma~\ref{Ysemicont}$(iii)$ gives $r>0$ with 
$$
\D(w;\widehat v)\le\tfrac12\D(u;\widehat v)=-\kappa_u\le0
$$ 
for all admissible $w$. Since $0\le\tfrac{\delta_0}{2}\|\widehat v-w\|_X$, the first inequality holds.

For the second inequality, suppose by contradiction no radius makes $\D(w;\widehat v)\le-\tfrac{\beta}{4}$ on the set $\left\{\F(w)\ge c-2\beta\right\}$. Then there exists $w_n\to u$ in $Y$, $\|w_n\|_X\le\Theta_0$, $\F(w_n)\ge c-2\beta$, with $\D(w_n;\widehat v)>-\tfrac{\beta}{4}$. Since $\phi(w_n)\to\phi(u)$ (see $(P_2)$),
$$
\liminf_{n\to\infty}\psi(w_n)=\liminf_{n\to\infty}(\F(w_n)-\phi(w_n))\ge(c-2\beta)-\phi(u),
$$
and $\langle\phi'(w_n),\widehat v-w_n\rangle\to\langle\phi'(u),\widehat v-u\rangle$. Hence
\begin{align*}
\limsup_{n\to\infty}\D(w_n;\widehat v) & \le\langle\phi'(u),\widehat v-u\rangle+\psi(\widehat v)-\liminf_n\psi(w_n)\\
	&\le\D(u;\widehat v)+\F(u)-(c-2\beta).
\end{align*}
Since $\F(u)<c-3\beta$ and $\D(u;\widehat v)\le0$, the right-hand side is smaller than $-\beta<-\tfrac{\beta}{4}$, so $\D(w_n;\widehat v)<-\tfrac{\beta}{4}$ for large $n$, a contradiction. Thus 
$$
\F(w)\ge c-2\beta\Rightarrow\D(w;\widehat v)\le-\tfrac{\beta}{4}=-\alpha.
$$

\item$(iii)$ Here $\F(u)<c-3\beta$ and $u\in\mathcal K_{\F}$. Then, there exists $\xi_u\in\partial\psi(u)$ with $\xi_u=-\phi'(u)$, so $\psi(u)-\psi(w)\le\langle\xi_u,u-w\rangle=\langle-\phi'(u),u-w\rangle$, whence
$$
\D(w;u)=\langle\phi'(w),u-w\rangle+\psi(u)-\psi(w)\le\langle\phi'(w)-\phi'(u),u-w\rangle .
$$
Property $(P_4)$ ensures the existence of $r>0$ such that $\|w-u\|_Y<r$ and $\|w\|_X\le\Theta_0$ imply $\|\phi'(w)-\phi'(u)\|_{X^\ast}\le\tfrac{\delta_0}{2}$. Choosing $\widehat v=u$,
$$
\D(w;\widehat v)=\D(w;u)\le\|\phi'(w)-\phi'(u)\|_{X^*}\|u-w\|_X\le\frac{\delta_0}{2}\|u-w\|_X=\frac{\delta_0}{2}\|\widehat v-w\|_X,
$$
so the first inequality holds. For the second, taking $\widehat v=u$ we have $\D(u;u)=0$, and the computation of case $(ii)$ (which used only $\F(u)<c-3\beta$ and $\D(u;\widehat v)\le0$) yields 
$$
\limsup_{n\to\infty}\D(w_n;u)\le\D(u;u)+\F(u)-(c-2\beta)<-\beta<-\tfrac{\beta}{2},
$$ so 
$$
\F(w)\ge c-2\beta\Rightarrow\D(w;\widehat v)\le-\tfrac{\beta}{2}\le-\alpha.
$$
\end{description}
In the cases $(ii)$ and $(iii)$ the bounds have been produced on two possibly different radii, so wechoose $r(u)$ to be the smaller of the two and thus both bounds hold simultaneously on the corresponding ball. Shrinking each $r(u)$ further we may assume $r(u)\le1$. In all three cases the stated bounds hold on $\left\{w:\|w\|_X\le\Theta_0,\ \|w-u\|_Y<r\right\}$.
\end{proof}

We can now assemble the deformation. The construction is run entirely with $Y$-balls, the compact embedding supplying the finite subcover that norm-compact linking sets would otherwise provide.

\begin{theorem}[Quantitative deformation for the Poincar\'{e} action]\label{QuantDef}
Assume condition $(C)$ holds and let $c\in\mathbb R$, $\beta\in(0,1)$, $\delta_0\in(0,\beta)$, $\rho>0$ and $Q\subseteq X$  (in the global case $Q=X$) be such that  $\slope{u}>5\delta_0$ for every $u\in D(\psi)$ with $|\F(u)-c|\le3\beta$ and $\dist_Y(u,Q)\le\rho$ $($equivalently, \eqref{SlopeStrip} holds$)$. Let
$$
K\subset\left\{u\in D(\psi):|\F(u)-c|\le3\beta,\ \dist_Y(u,Q)\le\rho/2\right\},
$$
with closure $\widehat K$ in $Y$, and let $\varrho_\ast\in\left(0,\min\{1,\rho/2\}\right]$. Then there exist $\varepsilon_\ast\in(0,\beta)$, a bounded open set ${\cal U}\subset X \cap [Q]^Y_{\rho}$ with $\widehat K\subset {\cal U}$, and maps $\widehat\sigma,\widetilde\sigma\in C([0,1]\times X,X)$ such that:
\begin{enumerate}[$(i)$]
\item $\widehat\sigma(0,\cdot)=\widetilde\sigma(0,\cdot)=\mathrm{id}$;
\item for every $t\in[0,1]$, $\widehat\sigma(t,\cdot)$ and $\widetilde\sigma(t,\cdot)$ are homeomorphisms of $X$ onto $X$ with jointly continuous inverses;
\item $\widehat\sigma(t,u)=\widetilde\sigma(t,u)=u$, whenever $u\notin {\cal U}$;
\item $\|\widehat\sigma(t,u)-u\|_X\le\varrho_\ast$ and $\|\widetilde\sigma(t,u)-u\|_X\le\varrho_\ast$, for all $(t,u)\in [0,1]\times X$.In particular, both maps displace every point by at most $\rho/2$ in the $Y$-metric;
\item $\F(\widehat\sigma(t,u))\le\max\{\F(u),c-\beta\}$, and
$$
\F(\widehat\sigma(t,u))=\F(u)\in(c-\beta,+\infty)\Rightarrow\widehat\sigma(t,u)=u;
$$
\item if $u\in K$ and $\F(u)\le c+\varepsilon_\ast$, then $\F(\widehat\sigma(1,u))\le c-\varepsilon_\ast$;
\item[$(\widetilde v)$] if $\F(u)\ge c-2\beta$, then $\F(\widetilde\sigma(t,u))\ge\min\{\F(u),c+\beta\}$, and 
$$
\F(\widetilde\sigma(t,u))=\F(u)<c+\beta\Rightarrow\widetilde\sigma(t,u)=u;
$$
\item[$(\widetilde{vi})$] if $u\in K$ and $\F(u)\ge c-\varepsilon_\ast$, then $\F(\widetilde\sigma(1,u))\ge c+\varepsilon_\ast$.
\end{enumerate}
\end{theorem}

\begin{proof}
If $K=\varnothing$ all the assertions are trivial: it suffices to take ${\cal U}=\varnothing$, $\widehat\sigma=\widetilde\sigma=\mathrm{id}$ and $\varepsilon_\ast=\beta/2$, since then $(i)$--$(v)$ and $(\widetilde v)$ hold with equality and no displacement, while $(vi)$ and $(\widetilde{vi})$ are vacuous. We may therefore assume $K\ne\varnothing$, hence $\widehat K\ne\varnothing$.

Let $\widehat R$ bound the strip $\Sigma_{c,\beta}$, so that $\widehat K\subset B_X\left[0,\widehat R\right]$, and put
$$
\Theta=\widehat R+2 .
$$
The proof proceeds in several steps as follows.

\begin{description}
\item{\sc Step} 1: {\it Selections and a finite cover.}

Applying Lemma~\ref{DescSel} with $\Theta_0=\Theta$, for each $u\in\widehat K$ it produces $\widehat v(u)\in D(\psi)$, $r(u)\in(0,1]$ and $\alpha(u)>0$.

We claim there exists $\delta>0$ such that for every $u\in\widehat K$ there exists $z\in\widehat K$ with the property that $B_Y(u,\delta)\subset B_Y\big(z,r(z)\big)$.

Arguing by contradiction, suppose this is not the case. Then there exists  a sequence $\{u_n\}\subset\widehat K$ such that $B_Y(u_n,1/n)$ is contained in no member of $\{B_Y(z,r(z))\}_{z\in\widehat K}$. As $\widehat K$ is compact in $Y$, a subsequence satisfies $\|u_{j_n}-u\|_Y\to0$ for some $u\in\widehat K$. Pick $n$ with $\|u_{j_n}-u\|_Y<\epsilon=r(u)/2$ and $1/j_n<\epsilon$. Then $B_Y(u_{j_n},1/j_n)\subset B_Y(u,2\epsilon)$, which is a member of the covering, since $u\in\widehat K$. This contradiction proves the claim.

Set
$$
\widehat\rho=\dist_Y\big(\widehat K,\,X\setminus[Q]^Y_{\rho}\big)\ \in\ [\rho/2,+\infty] ,
$$
with the convention $\dist_Y(\cdot,\varnothing)=+\infty$. If $u\in\widehat K$ and $z\notin[Q]^Y_{\rho}$, then
$$
\|u-z\|_Y\ \ge\ \dist_Y(z,Q)-\dist_Y(u,Q)\ >\ \rho-\rho/2\ =\ \rho/2 .
$$
In the global case $Q=X$ one has $X\setminus[Q]^Y_{\rho}=\varnothing$, hence $\widehat\rho=+\infty$. Define
\begin{equation}\label{Radius}
R=\min\big\{\delta,\ \tfrac13 \widehat\rho,\ 1\big\}\ \in\ (0,1] ,
\end{equation}
and note that, $R$ being smaller than $\delta$, the above claim holds with $R$ in place of $\delta$.

The $Y$-balls $\{B_Y(u,R/2)\}_{u\in\widehat K}$ cover the $Y$-compact set $\widehat K$, so we may extract a finite subcover, i.e., there exist $\widehat u_1,\dots,\widehat u_{N}\in\widehat K$ with
\begin{equation}\label{HalfCover}
\widehat K \subset \bigcup_{i=1}^N B_Y\left(\widehat u_i,R/2\right).
\end{equation}
For each $i\in\{1,\ldots,N\}$ choose a point $z_i\in\widehat K$ such that $B_Y(\widehat u_i,R)\subset B_Y\big(z_i,r(z_i)\big)$, and put
$$
v_i=\widehat v(z_i)\in D(\psi),\qquad \alpha_i=\alpha(z_i)>0 .
$$
Lemma~\ref{DescSel}, applied at the point $z_i$, then yields its two bounds we shall use in the sequel: for every $w\in X$ with $\|w\|\le\Theta$ and $\|w-\widehat u_i\|_Y<R$,
\begin{equation}\label{NeutralUnif}
\D(w;v_i)\ \le\ \tfrac{\delta_0}{2}\|v_i-w\|_X,
\end{equation}
and
\begin{equation}\label{HighUnif}
\F(w)\ge c-2\beta\ \Longrightarrow\ \D(w;v_i)\ \le\ -\alpha_i.
\end{equation}
Two consequences of \eqref{Radius} and \eqref{HalfCover} will be used repeatedly. First,
\begin{equation}\label{TubeCover}
X\cap\bigcup_{i=1}^N B_Y(\widehat u_i,R)\ \subset\ [Q]^Y_{\rho}.
\end{equation}
Indeed, if some $w\in X\cap B_Y(\widehat u_i,R)$ satisfied $w\notin[Q]^Y_{\rho}$, then $w$ would belong to $X\setminus[Q]^Y_{\rho}$ and $\widehat u_i$ to $\widehat K$, so that we would give $\widehat\rho\le\|w-\widehat u_i\|_Y<R\le\tfrac13 \widehat\rho$, which is impossible as $\widehat\rho>0$. For $Q=X$ the inclusion \eqref{TubeCover} is trivially true. Second, the set $[\widehat K]^Y_{R/2}$ is still covered, i.e.,
\begin{equation}\label{Fattened}
[\widehat K]^Y_{R/2}\ \subset\ \bigcup_{i=1}^N B_Y(\widehat u_i,R).
\end{equation}
Indeed, if $\dist_Y(w,\widehat K)\le R/2$, choose $u\in\widehat K$ with $\|w-u\|_Y\le R/2$  and, by \eqref{HalfCover}, an index $i$ with $\|u-\widehat u_i\|_Y<R/2$, which gives $\|w-\widehat u_i\|_Y<R$. Finally set
$$
\alpha=\min_{1\le i\le N}\alpha_i>0,\  D_0=\Theta+\max_{1\le i\le N}\|v_i\|_X,
$$
so that $\|v_i-w\|_X\le D_0$ whenever $\|w\|_X\le\Theta$ (note also $R\le1\le D_0$).

\item{\sc Step} 2: {\it Lipschitz partition of unity.}

For $i\in\{1,\ldots,N\}$ put
$$
p_i(w)=\begin{cases}
R-\|w-\widehat u_i\|_Y, & \mbox{if }\|w-\widehat u_i\|_Y<R,\\[2pt]
0, & \mbox{ otherwise,}
\end{cases}
$$
which is $1$-Lipschitz on $X$ (being $1$-Lipschitz in $\|\cdot\|_Y$). Define the radial cutoff
$$
\chi(w)=\begin{cases}
1, & \mbox{if }\|w\|_X\le\widehat R+\varrho_\ast,\\[2pt]
\widehat R+1+\varrho_\ast-\|w\|_X, & \mbox{if }\widehat R+\varrho_\ast<\|w\|_X<\widehat R+1+\varrho_\ast,\\[2pt]
0, & \mbox{if }\|w\|_X\ge\widehat R+1+\varrho_\ast,
\end{cases}
$$
which is also $1$-Lipschitz.

By \eqref{Fattened}, every point of $[\widehat K]^Y_{R/4}$ lies in one of the balls $B_Y(\widehat u_i,R)$, so that $\sum_{i=1}^N p_i>0$ there. In fact one has the explicit lower bound
\begin{equation}\label{MassLower}
\sum_{i=1}^N p_i(w)\ >\ \frac R4,\mbox{ for every }w\in[\widehat K]^Y_{R/4} .
\end{equation}
Indeed, given such a $w$, choose $u\in\widehat K$ with $\|w-u\|_Y\le R/4$ and, by \eqref{HalfCover}, an index $i$ with $\|u-\widehat u_i\|_Y<R/2$. Then $\|w-\widehat u_i\|_Y<3R/4$, whence $p_i(w)=R-\|w-\widehat u_i\|_Y>R/4$.
Put
$$
p_0(w)=\begin{cases}
\tfrac R4-\sum_{i=1}^N\chi(w)\,p_i(w), & \mbox{if }\sum_{i=1}^N\chi(w)\,p_i(w)<\tfrac R4,\\[2pt]
0, & \mbox{otherwise, }
\end{cases}
$$
and
$$
h_i(w)=\frac{\chi(w)\,p_i(w)}{p_0(w)+\sum_{j=1}^N\chi(w)\,p_j(w)}\ (i\in\{1,\ldots,N\}),\ 
m(w)=\sum_{i=1}^Nh_i(w) .
$$
The denominator is greater or equal to $R/4>0$ everywhere, so each $h_i$ is Lipschitz and $0\le m\le1$. Moreover:
\begin{description}
\item $\bullet$ $h_i(w)=0$ unless $\|w-\widehat u_i\|_Y<R$ and $\|w\|_X<\widehat R+1+\varrho_\ast$;
\item $\bullet$ $m\equiv1$ on $[\widehat K]^Y_{R/4}\cap\left\{\|w\|_X\le\widehat R+\varrho_\ast\right\}$: there $\chi=1$ and $\sum_{i=1}^N p_i>R/4$ by \eqref{MassLower}, hence $p_0=0$.
\end{description}

\item{\sc Step} 3: {\it The vector field and the elementary maps.}

Set
$$
{\cal T}(w)=\sum_{i=1}^Nh_i(w)\,\frac{v_i-w}{\max\left\{\|v_i-w\|_X,R\right\}} .
$$
Then ${\cal T}$ is globally Lipschitz on $X$ (a finite sum of Lipschitz terms, vanishing off a bounded set), $\|{\cal T}(w)\|_X\le m(w)\le1$, and ${\cal T}\equiv0$ off the bounded open set
$$
{\cal U}=\Big(X\cap\bigcup_{i=1}^N B_Y(\widehat u_i,R)\Big)\cap\left\{w\in X:\|w\|_X<\widehat R+1+\varrho_\ast\right\}\ \supset\ \widehat K.
$$
It is readily seen that $\widehat K\subset\bigcup_{i=1}^N B_Y(\widehat u_i,R)$ by \eqref{HalfCover}, $\widehat K\subset X$ and $\widehat K\subset\left\{\|w\|_X\le\widehat R\right\}$. Moreover ${\cal U}\subset[Q]^Y_{\rho}$ by \eqref{TubeCover}, since ${\cal U}\subset X\cap\bigcup_{i=1}^N B_Y(\widehat u_i,R)$ (a vacuous inclusion when $Q=X$).

Let $\nu_0$ be a Lipschitz constant of ${\cal T}$, and fix $\bar s\in(0,1]$ sufficiently small such that
$$
\bar s\nu_0\le\tfrac12,\quad \bar s\le\tfrac{\varrho_\ast}{2},\quad \bar s\le\tfrac{R}{2},\quad \bar s\le r_{\varepsilon_0}(\Theta),
$$
where $\varepsilon_0=\min\left\{\tfrac{\delta_0}{2},\tfrac{\alpha}{4D_0}\right\}$ and $r_{\varepsilon_0}(\Theta)$ is the radius from Lemma~\ref{ConvComb}. Then, for any $s\in [0,\bar s]$, the maps
$$
\eta_s=\mathrm{id}+s{\cal T},\  \bar\eta_s=\mathrm{id}-s{\cal T}
$$
are bi-Lipschitz homeomorphisms of $X$ with jointly continuous inverses, equal to the identity off ${\cal U}$, and $\|\eta_s(w)-w\|_X,\ \|\bar\eta_s(w)-w\|_X\le s$.

Let us prove these assertions. It is enough to treat $\eta_s$, since $\bar\eta_s$ is obtained by replacing ${\cal T}$ with $-{\cal T}$, which has the same Lipschitz constant $\nu_0$ and satisfies the same bound $\|-{\cal T}\|\le1$. Fix $s\in[0,\bar s]$.

For all $w_1,w_2\in X$, the Lipschitz continuity of ${\cal T}$ and $s\nu_0\le\bar s\nu_0\le\frac12$ give
\begin{equation}\label{BiLip}
\frac{1}{2}\|w_1-w_2\|_X \le \|\eta_s(w_1)-\eta_s(w_2)\|_X\le (1+s\nu_0)\|w_1-w_2\|_X\le\ \frac{3}{2}\|w_1-w_2\|_X,
\end{equation}
the middle inequality following from 
$$
\eta_s(w_1)-\eta_s(w_2)=(w_1-w_2)+s\big({\cal T}(w_1)-{\cal T}(w_2)\big)
$$ 
and $\|{\cal T}(w_1)-{\cal T}(w_2)\|_X\le\nu_0\|w_1-w_2\|_X$. In particular $\eta_s$ is injective.

Let $y\in X$ and consider $\Lambda_y:X\to X$, $\Lambda_y(w)=y-s{\cal T}(w)$. As 
$$
\|\Lambda_y(w_1)-\Lambda_y(w_2)\|_X\le s\nu_0\|w_1-w_2\|_X\le\frac{1}{2}\|w_1-w_2\|_X
$$ for all $w_1,w_2\in X$, the map $\Lambda_y$ is a contraction of the Banach space $X$, so it has a (unique) fixed point $w_0$ by the Banach fixed point theorem. Moreover, $w_0=y-s{\cal T}(w_0)$ means exactly $\eta_s(w_0)=y$. Together with \eqref{BiLip}, this shows that $\eta_s$ is a bi-Lipschitz homeomorphism of $X$ onto $X$, with $\eta_s^{-1}$ being $2$-Lipschitz on all of $X$.

Let $s_1,s_2\in[0,\bar s]$ and $y_1,y_2\in X$ be fixed and set $w_1=\eta_{s_1}^{-1}(y_1)$, $w_2=\eta_{s_2}^{-1}(y_2)$, so that $w_1=y_1-s_1{\cal T}(w_1)$ and $w_2=y_2-s_2{\cal T}(w_2)$. Subtracting,
$$
w_1-w_2=(y_1-y_2)-(s_1-s_2){\cal T}(w_1)-s_2\big({\cal T}(w_1)-{\cal T}(w_2)\big) ,
$$
hence, using $\|{\cal T}\|\le1$ and $s_2\nu_0\le\frac12$,
$$
\|w_1-w_2\|_X \le \|y_1-y_2\|_X+|s_1-s_2|+\frac{1}{2}\|w_1-w_2\|_X,
$$
that is,
$$
\big\|\eta_{s_1}^{-1}(y_1)-\eta_{s_2}^{-1}(y_2)\big\|_X \le 2\big(\|y_1-y_2\|_X+|s_1-s_2|\big) .
$$
Thus $(s,y)\mapsto\eta_s^{-1}(y)$ is Lipschitz, in particular continuous, on $[0,\bar s]\times X$.

\smallskip
Note that, if $w\notin {\cal U}$ then ${\cal T}(w)=0$, so $\eta_s(w)=\bar\eta_s(w)=w$. Thus, 
$$
\|\eta_s(w)-w\|_X=s\|{\cal T}(w)\|_X\le s
$$ because $\|{\cal T}(w)\|_X\le1$, and likewise for $\bar\eta_s$.

\item {\sc Step} 4: {\it One-step estimates.}

Throughout this step $s\in(0,\bar s]$ is an arbitrary parameter and $w\in D(\psi)$.

\smallskip
\noindent$(4a)$ \textit{The elementary step as a convex combination.}

Writing, 
$$
\theta_j=\frac{s h_j(w)}{\max\left\{\|v_j-w\|_X,R\right\}} \ge 0 ,
$$
the definition of ${\cal T}$ reads $s{\cal T}(w)=\sum_{j=1}^N\theta_j(v_j-w)$, so that
$$
\eta_s(w)=w+\sum_{j=1}^N\theta_j(v_j-w)=\left(1-{\sum_{j=1}^N}\theta_j\right)w+{\sum_{j=1}^N}\theta_j v_j 
$$
and 
$$
\bar\eta_s(w)=\Big(1+\sum_{j=1}^N\theta_j\Big)w-\sum_{j=1}^N\theta_j v_j .
$$
Then $\lambda=\sum_{j=1}^N \theta_j\leq \frac{1}{2}$ at every point $w\in X$. Indeed $\max\left\{\|v_j-w\|_X,R\right\}\ge R$ and $\sum_{j=1}^Nh_j(w)=m(w)\le1$, so
$$
\sum_{j=1}^N\theta_j \le\frac{s}{R}\sum_{j=1}^Nh_j(w) = \frac{s m(w)}{R}\le \frac{s}{R}\le \frac{1}{2},
$$
the last inequality being the requirement $s\le\bar s\le R/2$ imposed in Step~3. So $\eta_s(w)$ is a convex combination of the points $w,v_j\in D(\psi)$, and hence lies in $D(\psi)$, whenever $w\in D(\psi)$.

\smallskip
\noindent$(4b)$ {\it Interim energy estimates.} 

Assume next $m(w)>0$. Then $h_j(w)>0$ for at least one $j$, so $w\in {\cal U}$ by the definition of ${\cal U}$ in Step~3, and therefore
$$
\|w\|_X < \widehat R+1+\varrho_\ast\le\Theta ,
$$
the last inequality because $\varrho_\ast\le1$. This is precisely the localisation under which Lemma~\ref{ConvComb} and the two bounds \eqref{NeutralUnif}--\eqref{HighUnif} are available at $w$.

We may now invoke Lemma~\ref{ConvComb} with the point $w$, the vectors $v_j$ and the coefficients $\theta_j$, whose sum is smaller than or equal to $\frac{1}{2}$. Its smallness hypothesis is met, because $\|v_j-w\|_X\le\max\left\{\|v_j-w\|_X,R\right\}$ gives $\theta_j\|v_j-w\|_X\le sh_j(w)$ and therefore
$$
\sum_{j=1}^N\theta_j\|v_j-w\|_X\le s\sum_{j=1}^Nh_j(w)=s m(w)\le s \le r_{\varepsilon_0}(\Theta) .
$$
With $\varepsilon=\varepsilon_0$, estimate \eqref{R-Estimate} of Lemma~\ref{ConvComb} gives $|\mathcal{R}^\pm|\le\varepsilon_0\sum_{j=1}^N\theta_j\|v_j-w\|_X\le\varepsilon_0sm(w)$, so that \eqref{Forward} and \eqref{Backward} become
\begin{align}
\F(\eta_s(w))-\F(w)&\le\sum_{j=1}^N\theta_j\,\D(w;v_j)+\varepsilon_0sm(w),\label{Ledger}\\
\F(\bar\eta_s(w))-\F(w)&\ge-\sum_{j=1}^N\theta_j\,\D(w;v_j)-\varepsilon_0sm(w).\label{LedgerB}
\end{align}
It remains to estimate $\sum_{j=1}^N\theta_j\D(w;v_j)$, which is done in the next substeps.  Both use the following observation: if $\theta_j>0$ then $h_j(w)>0$, so the first bullet of Step~2 gives $\|w-\widehat u_j\|_Y<R$, and we have just seen that $\|w\|_X\le\Theta$. Thus \eqref{NeutralUnif} and \eqref{HighUnif} may be applied to every {\it active} index $j$.

\smallskip
\noindent$(4c)$ {\it The global bound.}

By \eqref{NeutralUnif}, $\D(w;v_j)\le\tfrac{\delta_0}{2}\|v_j-w\|_X$ for every active $j$. Multiplying by $\theta_j$ and summing, and using once more $\theta_j\|v_j-w\|_X\le sh_j(w)$,
$$
\sum_{j=1}^N\theta_j\,\D(w;v_j) \le \frac{\delta_0}{2}\sum_{j=1}^N\theta_j\|v_j-w\|_X\le\frac{\delta_0}{2}s m(w) \le \tfrac{\delta_0}{2} s ,
$$
so that \eqref{Ledger} and $\varepsilon_0\le\delta_0/2$ yield the following {\it global bound}
\begin{equation}\label{Global}
\F(\eta_s(w))-\F(w) \le\frac{\delta_0}{2}s+\varepsilon_0 s \le \delta_0s.
\end{equation}
Symmetrically, \eqref{LedgerB} gives $\F(\bar\eta_s(w))-\F(w)\ge-\delta_0 s$. If $m(w)=0$, then all the coefficients $h_i(w)$ vanish, hence ${\cal T}(w)=0$ and $\eta_s(w)=\bar\eta_s(w)=w$, so that both inequalities are trivially true. Hence \eqref{Global} and its counterpart hold at \textit{every} $w\in D(\psi)$.

\smallskip
\noindent$(4d)$ {\it The high-energy bounds.} 

Suppose in addition $\F(w)\ge c-2\beta$. Then for each  active index the high-energy bound \eqref{HighUnif}, $\D(w;v_j)\le-\alpha_j\le-\alpha$ holds. To exploit this we need a lower bound for the $\theta_j$ and it comes from $\max\left\{\|v_j-w\|_X,R\right\}\le D_0$, which holds because $\|w\|_X\le\Theta$ implies $\|v_j-w\|_X\le D_0$ because $R\le1\le D_0$. Thus $\theta_j\ge\frac{s}{D_0}h_j(w)$, whence $\sum_{j=1}^N\theta_j\ge\tfrac{s}{D_0}m(w)$ and
$$
\sum_{j=1}^N\theta_j\D(w;v_j) \le -\alpha\sum_{j=1}^N\theta_j \le -\frac{\alpha}{D_0}sm(w) .
$$
Inserting this in \eqref{Ledger} and using $\varepsilon_0\le\tfrac{\alpha}{4D_0}$,
$$
\F(\eta_s(w))-\F(w) \le -\tfrac{\alpha}{D_0}sm(w)+\tfrac{\alpha}{4D_0}sm(w)=-\tfrac{3\alpha}{4D_0}sm(w) ,
$$
and symmetrically from \eqref{LedgerB}. Setting $\delta_2=\tfrac{3\alpha}{8D_0}>0$, so that $\tfrac{3\alpha}{4D_0}=2\delta_2$, we obtain the \textit{high-energy bounds}
\begin{align}
\F(\eta_s(w))-\F(w)&\le-2\delta_2sm(w),\label{High}\\
\F(\bar\eta_s(w))-\F(w)&\ge\ 2\delta_2sm(w),\label{Back}
\end{align}
valid whenever $\F(w)\ge c-2\beta$ (again trivially so when $m(w)=0$: then ${\cal T}(w)=0$, so that $\eta_s(w)=\bar\eta_s(w)=w$ and \textit{both} sides of \eqref{High}, and of \eqref{Back}, vanish). Note that \eqref{High} produces a \textit{strict} decrease exactly when $m(w)>0$, and that the decrease is largest, namely $2\delta_2 s$, where $m\equiv1$.

\smallskip
\noindent$(4e)$ {\it The vector field does not vanish in the high-energy regime.}

The following consequence of Lemma~\ref{ConvComb} is worth noticing, as it describes exactly what the case $m(w)>0$ of $(4b)$ leaves open: 

\quad {\it if $w\in D(\psi)$ satisfies $m(w)>0$ and $\F(w)\ge c-2\beta$, then ${\cal T}(w)\ne0$.} 

\noindent Suppose by contradiction that ${\cal T}(w)=0$. Then 
$$
\sum_{j=1}^N\theta_j(v_j-w)=s{\cal T}(w)=0,$$ 
so that $\eta_s(w)=w$  and, by the definition of $\mathcal{R}^+$ in Lemma~\ref{ConvComb},$$
\mathcal{R}^+=\phi(\eta_s(w))-\phi(w)-\Big\langle\phi'(w),\sum_{j=1}^N\theta_j(v_j-w)\Big\rangle=0 .
$$
Consequently \eqref{Forward} reads
$$
0\ =\ \F(\eta_s(w))-\F(w)\ \le\ \sum_{j=1}^N\theta_j\,\D(w;v_j) .
$$
On the other hand $m(w)>0$ gives $\sum_{j=1}^N\theta_j\ge\tfrac{s}{D_0}m(w)>0$, so the high-energy estimate of $(4d)$ yields $\sum_{j=1}^N\theta_j\D(w;v_j)\le-\alpha\sum_{j=1}^N\theta_j<0$, a contradiction. Thus, wherever the deformation is meant to lower the energy, the elementary step does move the point.

\smallskip
\noindent$(4f)$ \textit{The backward map outside $D(\psi)$.}

Two remarks are in order about points where the energy is $+\infty$. First, if $\bar\eta_s(w)\notin D(\psi)$ while $w\in D(\psi)$, then $\F(\bar\eta_s(w))=+\infty$ and \eqref{Back} holds trivially. Second, $\bar\eta_s$ maps $X\setminus D(\psi)$ into itself. Indeed, let $w\notin D(\psi)$ and suppose, by contradiction, that $\bar\eta_s(w)\in D(\psi)$. Putting $\lambda=\sum_{j=1}^N\theta_j\le\frac12$, the identity $\bar\eta_s(w)=(1+\lambda)w-\sum_{j=1}^N\theta_j v_j$ may be rewritten as
$$
w\ =\ \frac{1}{1+\lambda} \bar\eta_s(w)+\sum_{j=1}^N\frac{\theta_j}{1+\lambda}v_j ,
$$
a convex combination of points of the convex set $D(\psi)$, hence $w\in D(\psi)$, a contradiction. Consequently, once a backward orbit leaves $D(\psi)$ its energy remains $+\infty$.

\item {\sc Step} 5: {\it Constructing the deformations.}

From now on we fix the step size $s=\bar s$, the threshold of Step~3. Since the estimates of Step~4 were proved there for an arbitrary parameter in $(0,\bar s]$, they are available both for $s$ and for every $s'\in(0,s]$. Put
$$
M=\lfloor\varrho_\ast/s\rfloor\ \ (\ge1, \mbox{ as }s\le\varrho_\ast/2),\mbox{ so that }  \frac{\varrho_\ast}{2}\le Ms\le\varrho_\ast ,
$$
and
$$
\varepsilon_\ast=\tfrac14\min\Big\{\beta,\ \delta_2\varrho_\ast,\ \tfrac{\delta_2R}{2c_0}\Big\}\in(0,\beta) .
$$
Define the deformation by iterating $\eta_s$ over $M$ equal subintervals,
$$
\widehat\sigma(t,\cdot)=\eta_{(Mt-(\ell-1))s}\circ\underbrace{\eta_s\circ\cdots\circ\eta_s}_{\ell-1}\ \text{for }t\in\Big[\tfrac{\ell-1}{M},\tfrac{\ell}{M}\Big],
$$
and $\widetilde\sigma$ likewise with $\bar\eta_s$ in place of $\eta_s$. Let us check the first four assertions of the theorem.

\begin{description}
\item $\bullet$ $(i)$ For $t=0$ we are in the subinterval corresponding to $\ell=1$ and the partial parameter is $0$, so $\widehat\sigma(0,\cdot)=\eta_0=\mathrm{id}$, likewise $\widetilde\sigma(0,\cdot)=\mathrm{id}$.

\item $\bullet$ \textit{The definition is consistent and jointly continuous.} At a node $t=\ell/M$ the two candidate expressions are $\eta_s\circ\eta_s^{\ell-1}$ and $\eta_0\circ\eta_s^{\ell}$, and both equal $\eta_s^\ell$. On each subinterval, $(t,u)\mapsto\eta_{(Mt-(\ell-1))s}(u)=u+(Mt-(\ell-1))s\,{\cal T}(u)$ is jointly continuous in $(t,u)$, and so is its composition with the fixed continuous map $\eta_s^{\ell-1}$. Hence the map $\widehat\sigma\in C([0,1]\times X,X)$, and similarly for $\widetilde\sigma$.

\item $\bullet$ $(ii)$ On each subinterval $\widehat\sigma(t,\cdot)$ is the composition of at most $M$ maps of the form $\eta_{s'}$ with $s'\in[0,s]\subset[0,\bar s]$, each of which is a bi-Lipschitz homeomorphism of $X$ onto $X$ by the assertions proved in Step~3, hence so is $\widehat\sigma(t,\cdot)$. Its inverse is the composition, in reverse order, of the maps $\eta_{s'}^{-1}$. Since $(s',y)\mapsto\eta_{s'}^{-1}(y)$ is jointly continuous by the same assertions and $t\mapsto(Mt-(\ell-1))s$ is continuous, the inverse depends jointly continuously on $(t,u)$ as well.

\item $\bullet$ $(iii)$ If $u\notin {\cal U}$ then ${\cal T}(u)=0$, so $\eta_{s'}(u)=u$ for every $s'$. The point is therefore fixed by every factor and $\widehat\sigma(t,u)=u$. The same holds for $\widetilde\sigma$.

\item $\bullet$ $(iv)$ Each factor displaces a point by at most $s$, and there are at most $M$ of them, so $\|\widehat\sigma(t,u)-u\|_X\le Ms\le\varrho_\ast$, and likewise for $\widetilde\sigma$.
\end{description}

\item {\sc Step} 6: {\it Proof of $(v)$.}

The whole point is the following {\it one-step} statement, in which $s'\in[0,s]$ is arbitrary:
\begin{equation}\label{OneStepV}
\F\big(\eta_{s'}(w)\big)\ \le\ \max\{\F(w),\,c-\beta\}\quad\mbox{for every }w\in X .
\end{equation}
Let us prove it. For $s'=0$ there is nothing to prove, $\eta_0=\mathrm{id}$. We may therefore assume $s'\in(0,s]$, for which the estimates of Step~4 are available. If $w\notin D(\psi)$, then $\F(w)=+\infty$ and there is nothing to prove either. Let then $w\in D(\psi)$ and distinguish two cases.

If $\F(w)\ge c-2\beta$, then \eqref{High} gives $\F(\eta_{s'}(w))\le\F(w)$ (the right-hand side of \eqref{High} being $\le0$), so \eqref{OneStepV} holds. If $\F(w)<c-2\beta$, then \eqref{Global} gives
$$
\F(\eta_{s'}(w))\ \le\ \F(w)+\delta_0 s'\ <\ (c-2\beta)+\beta\ =\ c-\beta ,
$$
because $\delta_0 s'\le\delta_0\le\beta$ (recall $s'\le\bar s\le1$ and $\delta_0<\beta$), and \eqref{OneStepV} holds again.

An immediate consequence of \eqref{OneStepV} is the following:
\begin{equation}\label{MuLevel}
\mu\ge c-\beta \mbox{ and }\\F(w)\le\mu\Longrightarrow\F\big(\eta_{s'}(w)\big) \le \mu\quad(0\le s'\le s) ,
\end{equation}
since in that case $\F\big(\eta_{s'}(w)\big)\le\max\{\F(w),c-\beta\}\le\max\{\mu,c-\beta\}=\mu$.

Note that for $t=0$ we have nothing to prove. Fix now  $t\in(0,1]$ and $u\in D(\psi)$. The points $u\notin D(\psi)$ may be left aside, the first assertion of $(v)$ being trivial there, its right-hand side equalling $+\infty$, and the second one vacuous, its hypothesis $\F(u)\in(c-\beta,+\infty)$ forcing $\F(u)<+\infty$. 

The intervals $\big(\frac{\ell-1}{M},\frac{\ell}{M}\big]$, $\ell=1,\dots,M$, cover $(0,1]$, so there exists an integer $n$ with $1\le n\le M$ and
$$
t\in\Big(\frac{n-1}{M},\frac{n}{M}\Big] .
$$
Define the $n+1$ points
$$
w_k=\widehat\sigma\Big(\frac{k}{M},u\Big)\ (0\le k\le n-1),\quad w_n=\widehat\sigma(t,u).
$$
We distinguish the following cases 
\begin{description}

\item {\sc Case $n=1$}.

Then $w_0=u$ and $w_1=\widehat \sigma(t,u)=\eta_{Mts}(u)$. Using \eqref{OneStepV} with $w=u$ and $s'=Mts$ we get 
$$
\F(w_1)\le \max\{\F(u),c-\beta\}, 
$$
which is exactly the first statement of $(v)$.

For the second statement, if $\F(w_1)=\F(u)>c-\beta>c-2\beta$, then by \eqref{High} we have 
$$
0=\F(w_1)-\F(u)\leq -2\delta_2 s m(u),
$$
which gives $m(u)=0$, hence ${\cal T}(u)=0$. Consequently, 
$$
\widehat \sigma(t,u)=\eta_{Mts}(u)=u+Mts{\cal T}(u)=u.
$$

\item {\sc Case $n\ge 2$}.

Define 
$$
s_k=s\ (0\le k\le n-2),\quad s_{n-1}=\big(Mt-(n-1)\big)s ,
$$
so that $s_k\in(0,s]$ for every $k$. Reading off the definition of $\widehat\sigma$ in Step~5, these points are nothing but the successive positions of $u$ under the elementary maps:
\begin{equation}\label{Chain6}
w_0=u,\qquad w_{k+1}=\eta_{s_k}(w_k)\ \ (0\le k\le n-1), \mbox{ so, } w_n=\widehat\sigma(t,u) .
\end{equation}
Indeed, the formula of Step~5 gives
$$
\widehat\sigma\Big(\tfrac{k}{M},\cdot\Big)=\underbrace{\eta_s\circ\cdots\circ\eta_s}_{k} \ (0\le k\le n-1),
\quad
\widehat\sigma(t,\cdot)=\eta_{s_{n-1}}\circ\underbrace{\eta_s\circ\cdots\circ\eta_s}_{n-1} ,
$$
the first equality because for $k\ge1$ the time $k/M$ lies in the subinterval corresponding to $\ell=k$, where the partial parameter equals $s$ (for $k=0$ it is simply $(i)$), and the second because $t$ lies in the $n$-th subinterval. Note that $u\in D(\psi)$ and that each $\eta_{s'}$ maps $D(\psi)$ into itself by $(4a)$, hence every $w_k$ lies in $D(\psi)$ and every value $\F(w_k)$ is finite. Finally, read along the chain \eqref{Chain6}, the one-step inequality \eqref{OneStepV}, with $w=w_k$ and $s'=s_k$, becomes
\begin{equation}\label{OneStepK}
\F(w_{k+1})\le \max\{\F(w_k),\,c-\beta\}\quad(0\le k\le n-1) .
\end{equation}
We claim that
\begin{equation}\label{NoRise}
\F(w_k) \le \mu_0\quad(0\le k\le n),
\end{equation}
where $\mu_0 =\max\{\F(u),c-\beta\}$.
This is true for $k=0$, since $\F(w_0)=\F(u)\le\mu_0$ by \eqref{Chain6}. Assume it holds for some $k\le n-1$. As $\mu_0\ge c-\beta$, the estimate \eqref{MuLevel}, applied at the point $w_k$ with the parameter $s_k\in[0,s]$, gives
$$
\F(w_{k+1})\ =\ \F\big(\eta_{s_k}(w_k)\big)\ \le\ \mu_0.
$$
The claim follows by induction on $k$ and this proves the first statement of $(v)$.

We turn our attention to the second assertion of $(v)$. Assume now, in addition, that
$$
\F\big(\widehat\sigma(t,u)\big) = \F(u) \in (c-\beta,+\infty) ,
$$
and let us prove that $\widehat\sigma(t,u)=u$. Then 
$
\mu_0= \F(u) ,
$
and the induction just performed says  
$$
\F(w_k) \leq \F(u), (0\le k\le n-1).
$$
The hypothesis pins moreover its two endpoints to that same value:
\begin{equation}\label{EndValues}
\F(w_0)= \F(w_n) = \F(u) .
\end{equation}
We claim now that $\F(w_k)\ge\F(u)$ for $k=n,n-1,\dots,0$. This is true for $k=n$ by \eqref{EndValues}. If $\F(w_{k+1})\ge\F(u)$ for some $k\le n-1$, then \eqref{OneStepK} gives
$$
\F(u)\le\F(w_{k+1}) \le \max\{\F(w_k),c-\beta\} .
$$
The maximum on the right cannot be $c-\beta$, since this would force $\F(u)\le c-\beta$, hence it equals $\F(w_k)$, and $\F(w_k)\ge\F(u)$. The claim follows now by descending induction on $k$.

\smallskip
Putting this together with \eqref{NoRise}
\begin{equation}\label{ConstEnergy}
\F(w_k) = \F(u)>c-\beta>c-2\beta\quad(0\le k\le n) .
\end{equation}

Fix $k\in\{0,\dots,n-1\}$. By \eqref{ConstEnergy} we have $\F(w_k)\ge c-2\beta$, so the high-energy bound \eqref{High} is available at $w_k$ with the parameter $s_k$. Since $\F(w_{k+1})=\F(w_k)$ by \eqref{ConstEnergy}, it yields
$$
0 = \F(w_{k+1})-\F(w_k) \le -2\delta_2s_km(w_k) \le0 ,
$$
that is, $m(w_k)=0$. As $m=\sum_{i=1}^Nh_i$ with $h_i\ge0$, all the coefficients $h_i(w_k)$ vanish, so ${\cal T}(w_k)=0$ and therefore $w_{k+1}=w_k+s_k{\cal T}(w_k)=w_k$ for all $0\le k\le n-1$. Consequently,
$$
\widehat\sigma(t,u) = w_n =w_{n-1} =\cdots=w_0 =u ,
$$
which is exactly the second statement of $(v)$.
\end{description}

\item {\sc Step} 7: {\it Proof of $(vi)$.}

Let $u\in K$, so $\|u\|_X\le\widehat R$, and $\F(u)\le c+\varepsilon_\ast$. Consider the full orbit 
$$
w_0=u, w_{\ell+1}=\eta_s(w_\ell), \ (0\le \ell\le M-1)
$$ and record two facts, valid for every $\ell\le M$.

\begin{description}
\item $\bullet$ {\it The orbit stays in the region where the cutoff $\chi$ equals $1$.} 

By property $(iv)$, every $w_\ell$ satisfies $\|w_\ell\|_X\le\widehat R+\varrho_\ast$.

\item $\bullet$ {\it If $w_\ell\in [\widehat K]^Y_{R/4}$, then $m\equiv 1$.} 

This follows from the previous point and the second bullet of Step~2,
so that a single step starting at such a $w_\ell$ decreases $\F$ by at least $2\delta_2s$, provided $\F(w_\ell)\ge c-2\beta$, by \eqref{High}.
\end{description}

We now distinguish three mutually exclusive cases.
\begin{enumerate}[$(a)$]

\item  {\it The energy drops below $c-2\beta$ at some time.} 

If $\F(w_{\ell_0})<c-2\beta$ for some $\ell_0\le M$, then successively applying \eqref{MuLevel}, with $s'=s$ and $\mu=c-\beta$, to the sub-orbit starting from $w_{\ell_0}$, gives
$$
\F(w_M) \le c-\beta \le c-\varepsilon_\ast ,
$$
because $\varepsilon_\ast\le\beta/4<\beta$.

\item {\it The energy never drops below $c-2\beta$, but the orbit leaves $[\widehat K]^Y_{R/4}$.} 

Assume $\F(w_\ell)\ge c-2\beta$ for all $\ell\le M$ and let $\ell_0\le M$ be the first index such that $w_{\ell_0}\notin[\widehat K]^Y_{R/4}$. Each step displaces a point by at most $s$ in $X$ (see end of Step 3), hence by at most $s$ in the $Y$-metric, so that
$$
\dist_Y(w_{\ell_0},\widehat K)\le \|w_{\ell_0}-u\|_Y\le\ell_0 s ,
$$
because $u\in K\subset\widehat K$. As $\dist_Y(w_{\ell_0},\widehat K)>R/4$, this forces
$$
\ell_0\ >\ \frac{R}{4s} .
$$
On the other hand, we know that 
$$
\F(w_\ell)\in [\widehat K]_{R/4}^Y, \mbox{ for }0\le\ell<\ell_0 
$$ 
and by the second bullet above, each of them decreases $\F$ by at least $2\delta_2 s$. Summing,
$$
\F(u)-\F(w_{\ell_0}) \ge2\delta_2 s\ell_0 >2\delta_2s\cdot\frac{R}{4s}= \frac{\delta_2R}{2} \ge 4\varepsilon_\ast ,
$$
the last inequality being the definition of $\varepsilon_\ast$. Hence 
$$
\F(w_{\ell_0})\le\F(u)-4\varepsilon_\ast\le c+\varepsilon_\ast-4\varepsilon_\ast\le c-\varepsilon_\ast.
$$ 
As in the previous case, applying \eqref{MuLevel}, with $s'=s$ and $\mu=c-\varepsilon_\ast\ge c-\beta$, to the sub-orbit starting from $w_{\ell_0}$ we get $\F(w_M)\le c-\varepsilon_\ast$.

\item {\it The energy never drops below $c-2\beta$ and the orbit stays in $[\widehat K]^Y_{R/4}$.} 

Then $m(w_\ell)=1$ for every $\ell<M$, and \eqref{High} applies at each of the $M$ steps, giving a total decrease of at least
$$
2\delta_2 M s \ge2\delta_2\cdot\frac{\varrho_\ast}{2}= \delta_2\varrho_\ast \ge 4\varepsilon_\ast ,
$$
since $Ms\ge\varrho_\ast/2$. Therefore $\F(w_M)\le c+\varepsilon_\ast-4\varepsilon_\ast\le c-\varepsilon_\ast$.
\end{enumerate}

In every case $\F(\widehat\sigma(1,u))=\F(w_M)\le c-\varepsilon_\ast$, which is $(vi)$.

\item {\sc Step} 8: {\it Proof of  $(\widetilde v)$.}

Assume $\F(u)\ge c-2\beta$. As in the previous step, note that for $t=0$ we have nothing to prove.

The one-step statement dual to \eqref{OneStepV} is even simpler here:
\begin{equation}\label{OneStepVt}
\F(w)\ge c-2\beta\Longrightarrow\F\big(\bar\eta_{s'}(w)\big) \ge \F(w)\quad(0< s'\le s) .
\end{equation}
In order to prove this, note that, if $w\in D(\psi)$ this is \eqref{Back}, whose right-hand side is nonnegative. If $w\notin D(\psi)$, then $\bar\eta_{s'}(w)\notin D(\psi)$ by $(4f)$, so that both sides equal $+\infty$.

The chain of elementary steps is produced exactly as in Step~6: since $\widetilde\sigma$ is built from $\bar\eta_s$ precisely as $\widehat\sigma$ is built from $\eta_s$, the same reasoning as in  Step~5 applies here. Choosing the integer $n$ and the parameters $s_0,\dots,s_{n-1}\in(0,s]$ as in Step~6 and putting $\widetilde w_k=\widetilde\sigma\big(\frac{k}{M},u\big)$ for $0\le k\le n-1$ and $\widetilde w_n=\widetilde\sigma(t,u)$, we obtain
\begin{equation}\label{Chain8}
\widetilde w_0=u,\  \widetilde w_{k+1}=\bar\eta_{s_k}(\widetilde w_k) \ (0\le k\le n-1).
\end{equation}

Since $\F(\widetilde w_0)=\F(u)\ge c-2\beta$, an induction using \eqref{OneStepVt} shows that the energy is at least  $ c-2\beta$ along the whole orbit and is non-decreasing there. Consequently
$$
\F(\widetilde\sigma(t,u)) \ge \F(u) \ge \min\{\F(u),c+\beta\} ,
$$
which is the first half of $(\widetilde v)$.

We turn to the second assertion of $(\widetilde v)$. Fix $t\in(0,1]$ and $u\in X$ with
$$
\F\big(\widetilde\sigma(t,u)\big)= \F(u) \in [c-2\beta,\,c+\beta).
$$
In particular, $\F(u)<+\infty$, i.e., $u\in D(\psi)$. Let us prove that $\widetilde\sigma(t,u)=u$.

Since $\F(\widetilde w_n)=\F(\widetilde w_0)$, we claim the energy is constant along $\widetilde w_0,\widetilde w_1,\dots,\widetilde w_n$. 
Now $\F(\widetilde w_0)=\F(u)\ge c-2\beta$, and \eqref{OneStepVt}, applied at $\widetilde w_k$ with the parameter $s_k$, turns $\F(\widetilde w_k)\ge c-2\beta$ into $\F(\widetilde w_{k+1})\ge\F(\widetilde w_k)\ge c-2\beta$. Induction on $k$ therefore yields at once
$$
\F(\widetilde w_k)\ge c-2\beta \ (0\le k\le n)\mbox{ and } \F(\widetilde w_0)\le\F(\widetilde w_1)\le\cdots\le\F(\widetilde w_n) .
$$
By \eqref{Chain8} the first and the last term of this non-decreasing sequence are equal, hence it is constant:
\begin{equation}\label{ConstEnergy8}
\F(\widetilde w_k)= \F(u)\in [c-2\beta,\,c+\beta)\quad(0\le k\le n) .
\end{equation}
In particular these values are finite, so that every $\widetilde w_k$ lies in $D(\psi)$. This has to be checked here, whereas in Step~6 it was automatic, the forward map preserving $D(\psi)$ by $(4a)$ while $\bar\eta_s$ need not. 

Fix $k\in\{0,\dots,n-1\}$. Then by \eqref{ConstEnergy8} we have $\widetilde w_k\in D(\psi)$ and $\F(\widetilde w_k)\ge c-2\beta$, so \eqref{Back} is available at $\widetilde w_k$ with the parameter $s_k$. Since $\F(\widetilde w_{k+1})=\F(\widetilde w_k)$ by \eqref{ConstEnergy8}, it yields
$$
0\ =\ \F(\widetilde w_{k+1})-\F(\widetilde w_k)\ \ge\ 2\delta_2\,s_k\,m(\widetilde w_k)\ \ge\ 0 ,
$$
that is, $m(\widetilde w_k)=0$, whence ${\cal T}(\widetilde w_k)=0$ and $\widetilde w_{k+1}=\widetilde w_k-s_k{\cal T}(\widetilde w_k)=\widetilde w_k$ and \eqref{Chain8} gives
$$
\widetilde\sigma(t,u)=\widetilde w_n = \widetilde w_{n-1} =\cdots= \widetilde w_0 =u ,
$$
which is exactly the second statement of $(\widetilde v)$.

\item {\sc Step} 9: {\it Proof of $(\widetilde{vi})$.}

Let $u\in K$ with $\F(u)\ge c-\varepsilon_\ast\ (\ge c-2\beta)$, and consider the full orbit 
$$
\widetilde w_0=u,\  \widetilde w_{\ell+1}=\bar\eta_s(w_\ell), \ (0\le \ell \le M-1)\mbox{ so } \widetilde w_M=\widetilde\sigma(1,u).
$$
 As in Step~7, $\|\widetilde w_\ell\|_X\le\widehat R+\varrho_\ast$ for every $\ell$, so that $m(\widetilde w_\ell)=1$ whenever $\widetilde w_\ell\in[\widehat K]^Y_{R/4}$. By $(\widetilde v)$ the value never drops below $c-2\beta$, so \eqref{Back} applies at every step (an exit from $D(\psi)$ counting as a jump to $+\infty$). Two cases arise, and the situation is exactly the one of Step~7, with decreases replaced by increases.
\begin{enumerate}[$(a')$]

\item  \textit{The orbit leaves $[\widehat K]^Y_{R/4}$.} 

Let $\ell_0$ be the first index with $\widetilde w_{\ell_0}\notin[\widehat K]^Y_{R/4}$. As in case $(b)$ of Step~7 one gets $\ell_0>R/(4s)$, and each of the first $\ell_0$ steps increases $\F$ by at least $2\delta_2s$. Thus
$$
\F(\widetilde w_{\ell_0})-\F(u) \ge 2\delta_2s \ell_0 >\frac{\delta_2R}{2c_0} \ge 4\varepsilon_\ast ,
$$
so that $\F(\widetilde w_{\ell_0})\ge c-\varepsilon_\ast+4\varepsilon_\ast=c+3\varepsilon_\ast\ge c+\varepsilon_\ast$. Applying $(\widetilde v)$ to the sub-orbit starting at $\widetilde w_{\ell_0}$ then gives
$$
\F(\widetilde w_M) \ge\min\{\F(\widetilde w_{\ell_0}),c+\beta\} \ge \min\{c+\varepsilon_\ast,c+\beta\}=c+\varepsilon_\ast ,
$$
because $\varepsilon_\ast<\beta$.

\item \textit{The orbit stays in $[\widehat K]^Y_{R/4}$.} 

Then $m(\widetilde w_\ell)=1$ for every $\ell<M$ and \eqref{Back} gives a total increase of at least $2\delta_2Ms\ge\delta_2\varrho_\ast\ge4\varepsilon_\ast$, so that
$$
\F(\widetilde w_M) \ge \F(u)+4\varepsilon_\ast\ge c-\varepsilon_\ast+4\varepsilon_\ast \ge c+\varepsilon_\ast .
$$
\end{enumerate}
In either case $\F(\widetilde\sigma(1,u))\ge c+\varepsilon_\ast$, which is $(\widetilde{vi})$.
\end{description}
\end{proof}

\section{Proof of Minimax theorem for the Poincar\'{e} action functional}\label{SecMinimax}

We now prove Theorem \ref{MinMaxThm} using the deformation from the previous section. The main idea is that assuming $\F$ possesses no $(PS)$-sequence at the minimax levels $c$ and $d$, respectively, Lemma~\ref{Equiv} activates the deformation of Theorem~\ref{QuantDef}, which lowers (or raises) a near-optimal competitor below (above) the level, contradicting its minimax definition.

\begin{proof}[Proof of Theorem~\ref{MinMaxThm}]
We proceed in several steps as follows.

\begin{description}
\item {\sc Step} 1:{\it The four levels are finite and ordered. More precisely,}

\begin{equation}\label{Chain}
-\infty<d\le a\le b\le c<+\infty .
\end{equation}
Since $\psi\ge0$, $\F\ge\phi$ on $D(\psi)$ and $\F=+\infty$ off $D(\psi)$. By $(P_6)$, evaluating $\tilde\gamma$ at $t=0$ gives 
$$
\F(u)=\F(\tilde\gamma(0,u))\le\sup_{\substack{t\in [0,1]\\ z\in A}}\F(\tilde\gamma(t,z))<\infty, \mbox{ for every }u\in A.
$$ 
On the other hand, since $A$ is bounded and $\phi$ maps bounded sets to bounded sets 
$$
-\infty<\inf_A\phi\le a=\sup_A\F<\infty.
$$ 
Moreover, 
$$
c=\inf_{\gamma\in\HH}\sup_{\substack{t\in[0,1]\\u\in A}}\F(\gamma(t,u))\le\sup_{\substack{t\in[0,1]\\u\in A}}\F(\tilde\gamma(t,u))<\infty,
$$
by property $(P_6)$. By duality, $(P_6)$ gives 
$$
d=\sup_{\gamma\in\HH}\inf_{\substack{t\in[0,1]\\u\in B}}\F(\gamma(t,u))\ge\inf_{\substack{t\in[0,1]\\u\in B}}\F(\tilde\gamma(t,u))>-\infty.
$$

For the linking bounds fix $\gamma\in\HH$. Since $A$ links $B$ there exist $t_0\in(0,1]$, $u_0\in A$ with $\gamma(t_0,u_0)\in B$, so 
$$
\sup_{\substack{t\in[0,1]\\u\in A}}\F(\gamma(t,u))\ge\F(\gamma(t_0,u_0))\ge\inf_B\F=b.
$$
Taking infimum over $\gamma\in \HH$ gives $c\ge b$. In particular, $b\le c<\infty$. Since $B$ links $A$ there exist $t_1\in(0,1]$, $u_1\in B$ with $\gamma(t_1,u_1)\in A$, so 
$$
\inf_{\substack{t\in[0,1]\\u\in B}}\F(\gamma(t,u))\le\F(\gamma(t_1,u_1))\le\sup_A\F=a
$$
and taking supremum over $\gamma\in \HH$ gives $d\le a$. Combining with $a\le b$ from condition $(S)$ yields \eqref{Chain}.

\item {\scshape Step} 2:\textit{ Existence of the critical point $u$ with $\F(u)=c\ge b$.}

We distinguish the following cases.

\begin{description}

\item {\scshape Case 2A:} {$c>b$.}

 Then $a\le b<c$. Fix $\beta_0\in(0,1)$ with $a<c-3\beta_0$, and suppose, by contradiction, that $\F$ possesses no $(PS)_c$-sequence.

By Lemma~\ref{Equiv} (with $Q=X$) there are $\beta\in(0,1)$ and $\delta_0\in(0,\beta)$ for which \eqref{SlopeStrip} holds at level $c$ with $Q=X$ and an arbitrary $\rho>0$, say $\rho=2$. Recall that for $Q=X$ the constraint $\dist_Y(u,Q)\le\rho$ is vacuous, so the value of $\rho$ is immaterial here, while this choice gives $\min\{1,\rho/2\}=1$. Shrinking $\beta$ we may assume $\beta\le\beta_0$: a smaller strip inherits the same selection $v(u)$, and $\delta_0$ may be replaced by $\min\{\delta_0,\beta/2\}$. In particular, $a<c-3\beta$. Applying Theorem~\ref{QuantDef} with $Q=X$, $K=\Sigma_{c,\beta}$ and $\varrho_\ast=1$ produces $\varepsilon_\ast\in(0,\beta)$ and a deformation $\widehat\sigma$ satisfying $(i)$--$(vi)$.

Fix $\varepsilon\in (0,\varepsilon_\ast)$. Then by the definition of $c$, there exists $\gamma_\varepsilon\in \HH$ such that 
$$
\sup_{\substack{t\in[0,1]\\u\in A}}\F(\gamma_\varepsilon(t,u))\le c+\varepsilon.
$$
Define
$$
\gamma^*(s,u)=\begin{cases}\widehat\sigma(2s,u), & 0\le s\le\tfrac12,\\[2pt]\widehat\sigma\big(1,\gamma_\varepsilon(2s-1,u)\big), & \frac{1}{2}\le s\le1.\end{cases}
$$
Then $\gamma^*\in\HH$. Indeed $\gamma^*(0,\cdot)=\mathrm{id}$, for $s<1$ the mapping $\gamma^*(s,\cdot)$ is a homeomorphism ($\widehat\sigma(2s,\cdot)$ for $s<\tfrac12$, and the composition $\widehat\sigma(1,\cdot)\circ\gamma_\varepsilon(2s-1,\cdot)$ for $\tfrac12\le s<1$), with the two branches agreeing at $s=\tfrac12$ and with jointly continuous inverse. Since $\gamma_\varepsilon(1,\cdot)\equiv\zeta$ for some $\zeta\in X$, we get $\gamma^*(1,\cdot)\equiv\widehat\sigma(1,\zeta)$. 

Finally, the convergence $\gamma^*(s,u)\to\widehat\sigma(1,\zeta)$ as $s\to1$ is uniform on bounded sets: the continuity of $\widehat\sigma(1,\cdot)$ at the single point $\zeta$, guaranteed by $(ii)$, already suffices. Indeed, let $P\subset X$ be bounded and $\tau>0$. Then there exists $\delta>0$ such that 
$$
\|z-\zeta\|_X<\delta \Rightarrow\|\widehat\sigma(1,z)-\widehat\sigma(1,\zeta)\|_X<\tau,
$$ 
and, $\gamma_\varepsilon$ being a member of $\HH$, there exists $t_1\in[0,1)$ with $\|\gamma_\varepsilon(t,u)-\zeta\|_X<\delta$ for all $t\in(t_1,1]$ and $u\in P$. Hence $\|\gamma^*(s,u)-\widehat\sigma(1,\zeta)\|_X<\tau$ for all $s\in\big(\tfrac{1+t_1}{2},1\big]$ and $u\in P$.

We now bound $\F(\gamma^*(s,u))$ for $u\in A$.

If $s\le\frac{1}{2}$, then $\F(u)\le a<c-\beta$, so $(v)$ gives
$$
\F(\widehat\sigma(2s,u)) \le \max\{\F(u),c-\beta\} =c-\beta\le c-\varepsilon_\ast .
$$

If $s\ge\frac{1}{2}$, put $w=\gamma_\varepsilon(2s-1,u)$, so $\F(w)\le c+\varepsilon\le c+\varepsilon_\ast$. 

When $\F(w)\le c-\varepsilon_\ast$, property $(v)$ gives $\F(\widehat\sigma(1,w))\le c-\varepsilon_\ast$. If $\F(w)>c-\varepsilon_\ast$, then $|\F(w)-c|\le\varepsilon_\ast\le3\beta$, so $w\in K$ and $(vi)$ gives $\F(\widehat\sigma(1,w))\le c-\varepsilon_\ast$.

In all cases $\F(\gamma^*(s,u))\le c-\varepsilon_\ast$, whence
$$
c=\inf_{\gamma\in\HH} \sup_{\substack{t\in[0,1]\\u\in A}}\F(\gamma(t,u)) \le\sup_{\substack{s\in[0,1]\\u\in A}}\F(\gamma^*(s,u))\le c-\varepsilon_\ast ,
$$
a contradiction.

Therefore $\F$ possesses a $(PS)_c$-sequence $\{u_n\}$. By Lemma~\ref{PScritical} it is bounded and a subsequence converges in $Y$ to some $u\in\mathcal K_{\F}(c)$.

\smallskip
\item {\scshape Case 2B:} {$c=b$.}

 Here $c=\inf_B\F$, so $B$ realises the level and we seek a $(PS)_c$-sequence approaching $B$ in $Y$. Suppose none exists.

By Lemma~\ref{Equiv} with $Q=B$ there are $\beta\in(0,1)$, $\delta_0\in(0,\beta)$ and $\rho>0$ for which \eqref{SlopeStrip} holds at level $c$ with $Q=B$. Apply Theorem~\ref{QuantDef} with this $Q=B$, radius $\rho$, set
$$
K=\{u\in D(\psi):\ |\F(u)-c|\le3\beta\}\cap[B]^Y_{\rho/2},
$$
and $\varrho_\ast=\min\{1,\rho/c\}$, which is an admissible choice in Theorem~\ref{QuantDef}. This yields $\varepsilon_\ast\in(0,\beta)$ and a deformation $\widehat\sigma$ that is the identity off $[B]^Y_\rho$ and, by $(iv)$, moves points by at most $\rho/2$ in $Y$.

Fix $\varepsilon\in(0,\varepsilon_\ast)$ and take $\gamma_\varepsilon \in\HH$ with $\sup_{\substack{t\in[0,1]\\u\in A}}\F(\gamma_\varepsilon(t,u))\le c+\varepsilon$, then form $\gamma^*$ as in Case 2A. Since $A$ links $B$, there exists $(s^*,y)\in(0,1]\times A$ with $q=\gamma^*(s^*,y)\in B$.

Suppose first $s^*\ge\tfrac12$, and set $w=\gamma_\varepsilon(2s^*-1,y)$, so $q=\widehat\sigma(1,w)\in B$ and 
$$
\F(q)\ge\inf_B\F=b=c.
$$ 
By $(v)$, $c\le\F(q)\le\max\{\F(w),c-\beta\}$, which forces $\F(w)\ge c$, which combined with $\F(w)\le c+\varepsilon_\ast$ gives $|\F(w)-c|\le\varepsilon_\ast\le3\beta$. Moreover, $\widehat\sigma$ moves points by at most $\rho/2$ in $Y$ by $(iv)$, so
$$
\dist_Y(w,B)\ \le\ \dist_Y(q,B)+\rho/2\ =\ \rho/2 ,
$$
whence $w\in K$. Then $\F(w)\le c+\varepsilon_\ast$ and $(vi)$ give $\F(q)=\F(\widehat\sigma(1,w))\le c-\varepsilon_\ast<c$, contradicting $q\in B$.

Suppose next $s^*\le\tfrac12$, so $q=\widehat\sigma(2s^*,y)\in B$ with $y\in A$ and $\F(y)\le a\le b=c$. By $(v)$, $\F(q)\le\max\{\F(y),c-\beta\}\le c$, while $q\in B$ forces $\F(q)\ge c$, hence $\F(q)=c$. Thus $\max\{\F(y),c-\beta\}=c$, and, since $c-\beta<c$, forces $\F(y)=c>c-\beta$. Then the second assertion of $(v)$ gives $q=y\in A$, contradicting $A\cap B=\varnothing$.

Both cases are impossible, so a $(PS)_c$-sequence $\{u_n\}$ with $\dist_Y(u_n,B)\to0$ exists. By Lemma~\ref{PScritical} it is bounded and a subsequence converges in $Y$ to some $u\in{\mathcal K}_{\F}(c)$. Finally $\dist_Y(u_n,B)\to0$ yields $\dist_Y(u,B)=0$, i.e., $u\in\overline B^{\,Y}$.

In either case $\F(u)=c\ge b$, and $u\in\overline B^{\,Y}$ whenever $c=b$.
\end{description}

\item {\scshape Step} 3:\textit{Existence of the critical point $v$ with $\F(v)=d\le a$.} 

This is the mirror of Step 2 for the dual minimax level 
$$
d=\sup_{\gamma\in\HH}\inf_{\substack{t\in[0,1]\\u\in B}}\F(\gamma(t,u)),
$$
using the deformation $\widetilde\sigma$ of Theorem~\ref{QuantDef} (properties $(\widetilde v),(\widetilde{vi})$) and near-optimal paths $\gamma_\varepsilon\in \HH$ with $\inf_{\substack{t\in[0,1]\\u\in B}}\F(\gamma_\varepsilon(t,u))\ge d-\varepsilon$, the parameter $\varepsilon$ being fixed in each of the two cases once $\varepsilon_\ast$ is available. We indicate the two cases.

\begin{description}

\item {\scshape Case 3A:} $d<a$.

Then $d<b$. Fix $\beta_0\in(0,1)$ with $b>d+3\beta_0$, and suppose $\F$ possesses no $(PS)_d$-sequence. As in Case 2A, Lemma~\ref{Equiv} with $Q=X$ furnishes $\beta\in(0,1)$ with $\beta\le\beta_0$ and $\delta_0\in(0,\beta)$ with \eqref{SlopeStrip} at level $d$ (again with an arbitrary $\rho>0$, e.g., $\rho=2$, the constraint $\dist_Y(u,X)\le\rho$ being vacuous), so $b>d+3\beta$. Theorem~\ref{QuantDef} (with $Q=X$, $\rho=2$, $K=\Sigma_{d,\beta}=\{u\in D(\psi):|\F(u)-d|\le3\beta\}$ and $\varrho_\ast=1$) gives $\varepsilon_\ast\in(0,\beta)$ and the deformation $\widetilde\sigma$ which satisfies $(\widetilde v)$ and $(\widetilde{vi})$. Fix $\varepsilon\in(0,\varepsilon_\ast)$ and $\gamma_\varepsilon\in\HH$ with $\inf_{t\in[0,1],\,u\in B}\F(\gamma_\varepsilon(t,u))\ge d-\varepsilon$. Define
$$
\gamma^*(s,u)=\begin{cases}\widetilde\sigma(2s,u), & 0\le s\le\tfrac12,\\[2pt]\widetilde\sigma\big(1,\gamma_\varepsilon(2s-1,u)\big), & \tfrac12\le s\le1,\end{cases}\ \in\HH
$$
(that $\gamma^*\in\HH$ is checked exactly as in Case 2A, with $\widetilde\sigma$ in place of $\widehat\sigma$) and fix $u\in B$.

If $s\le\tfrac12$, then $\F(u)\ge b>d+3\beta>d-2\beta$, so $(\widetilde v)$ gives
$$
\F(\widetilde\sigma(2s,u)) \ge \min\{\F(u),d+\beta\}= d+\beta \ge d+\varepsilon_\ast .
$$

If $s\ge\tfrac12$, put $w=\gamma_\varepsilon(2s-1,u)$, so $\F(w)\ge d-\varepsilon\ge d-\varepsilon_\ast\ge d-2\beta$. 

If $\F(w)\ge d+\varepsilon_\ast$, then $(\widetilde v)$ gives $\F(\widetilde\sigma(1,w))\ge d+\varepsilon_\ast$, otherwise $|\F(w)-d|\le\varepsilon_\ast\le3\beta$, so $w\in K$ and $(\widetilde{vi})$ gives $\F(\widetilde\sigma(1,w))\ge d+\varepsilon_\ast$.

Hence $\inf_{\substack{s\in[0,1]\\u\in B}}\F(\gamma^*(s,u))\ge d+\varepsilon_\ast$, contradicting the definition of $d$. Therefore $\F$ possesses a $(PS)_d$-sequence, and Lemma~\ref{PScritical} yields $v\in\mathcal K_{\F}(d)$.

\smallskip
\item {\scshape Case 3B:} $d=a$.

Here $d=\sup_A\F$, so $A$ is the dual set. Suppose no $(PS)_d$-sequence with $\dist_Y(\cdot,A)\to0$ exists.

Lemma~\ref{Equiv} with $S=A$ gives $\beta\in(0,1)$, $\delta_0\in(0,\beta)$ and $\rho>0$ for which \eqref{SlopeStrip} holds at level $d$ with $Q=A$. Theorem~\ref{QuantDef} with $Q=A$, radius $\rho$,
$$
K=\{u\in D(\psi):\ |\F(u)-d|\le3\beta\}\cap[A]^Y_{\rho/2}
$$ and $\varrho_\ast=\min\{1,\rho/2\}$ gives $\varepsilon_\ast\in(0,\beta)$ and $\widetilde\sigma$ that is the identity off $[A]^Y_\rho$ and, by $(iv)$, moves points by at most $\rho/2$ in $Y$. Fix $\varepsilon\in(0,\varepsilon_\ast)$ and $\gamma_\varepsilon\in\HH$ with 
$$
\inf_{t\in[0,1],\,u\in B}\F(\gamma_\varepsilon(t,u))\ge d-\varepsilon.
$$ Set $\gamma^*$ as above. Since $B$ links $A$, there exists $(s^*,y)\in(0,1]\times B$ such that $q=\gamma^*(s^*,y)\in A$, so $\F(q)\le\sup_A\F=d$.

Suppose first $s^*\ge\tfrac12$, and set $w=\gamma_\varepsilon(2s^*-1,y)$, so $q=\widetilde\sigma(1,w)$ and 
$$
\F(w)\ge d-\varepsilon\ge d-\varepsilon_\ast\ge d-2\beta.
$$ By $(\widetilde v)$, $\F(q)\ge\min\{\F(w),d+\beta\}$. Since $\F(q)\le d$, we must have $\F(w)<d+\varepsilon_\ast$, otherwise $\F(q)\ge\min\{d+\varepsilon_\ast,d+\beta\}=d+\varepsilon_\ast>d$, a contradiction. Hence $|\F(w)-d|\le\varepsilon_\ast\le3\beta$. Moreover $\widetilde\sigma$ moves points by at most $\rho/2$ in $Y$ by $(iv)$, so
$$
\dist_Y(w,A)\ \le\ \dist_Y(q,A)+\rho/2\ =\ \rho/2 ,
$$
whence $w\in K$. Then $(\widetilde{vi})$ gives $\F(q)=\F(\widetilde\sigma(1,w))\ge d+\varepsilon_\ast>d$, contradicting the fact that $q\in A$.

Suppose next $s^*\le\tfrac12$, so $q=\widetilde\sigma(2s^*,y)\in A$ with $y\in B$ and $\F(y)\ge b\ge d\ge d-2\beta$. By $(\widetilde v)$, $\F(q)\ge\min\{\F(y),d+\beta\}\ge\min\{d,d+\beta\}=d$. But, $q\in A$ and this forces $\F(q)\le d$. Thus $\F(q)=d=\min\{\F(y),d+\beta\}$, so $\F(y)=d<d+\beta$. The second assertion of $(\widetilde v)$ then gives $q=y\in B$, contradicting $A\cap B=\varnothing$.

Both cases are impossible, so a $(PS)_d$-sequence with $\dist_Y(\cdot,A)\to0$ must exist. As in Case 2B, Lemma~\ref{PScritical} provides $v\in \mathcal K_{\F}(d)$ and $v\in\overline A^{\,Y}$.

\end{description}

In either case $\F(v)=d\le a$, and $v\in\overline A^{\,Y}$ whenever $d=a$.

\item {\scshape Step} 4:\textit{ The critical points are distinct, provided $\overline A^{\,Y}\cap\overline B^{\,Y}=\varnothing$.} 

We have $\F(u)=c$, $\F(v)=d$, and \eqref{Chain} gives $d\le a\le b\le c$. If $c>d$ then $\F(u)>\F(v)$, so $u\ne v$. If $c=d$, the chain collapses to $d=a=b=c$. In particular, $c=b$, so Case 2B gives $u\in\overline B^{\,Y}$, and $d=a$, so Case 3B gives $v\in\overline A^{\,Y}$. Thus, when $c=d$, $u\in\overline B^{\,Y}$ and $v\in\overline A^{\,Y}$. Consequently,  if  $\overline A^{\,Y}\cap\overline B^{\,Y}=\varnothing$, then these points lie in disjoint sets and $u\ne v$.
\end{description}
\end{proof}

\pagestyle{plain}
\bibliographystyle{siam}
\bibliography{biblio}

\end{document}